\documentclass[article,onesided]{amsart}\newif\ifprivate\privatefalse \newif\iflongversion\longversionfalse \newif\ifoldversion\oldversionfalse \def\loadTIKZ{\usepackage{tikz}\usetikzlibrary{matrix,arrows,calc,cd,decorations.pathmorphing}}\overfullrule=5pt \newcommand\hfuzzReset{\hfuzz=3pt}\hfuzzReset \newcommand\toleranceReset{\tolerance=1400}\toleranceReset \newcommand\emergencystretchReset{\emergencystretch=1ex}\emergencystretchReset \hbadness=10000 \usepackage{ifpdf}\newcommand{\emailTressl}{marcus.tressl@manchester.ac.uk}\newcommand{\homepageTressl}{\url{http://personalpages.manchester.ac.uk/staff/Marcus.Tressl/}}\let\oldtocsection=\tocsection \let\oldtocsubsection=\tocsubsection \let\oldtocsubsubsection=\tocsubsubsection \renewcommand{\tocsection}[2]{\hspace{0em}\vspace*{0.1em}\oldtocsection{#1}{#2}}\renewcommand{\tocsubsection}[2]{\hspace{4ex}\oldtocsubsection{#1}{#2}}\renewcommand{\tocsubsubsection}[2]{\hspace{6ex}\oldtocsubsubsection{#1}{#2}}\ifpdf \usepackage[pdftex]{lscape}\else \usepackage{lscape}\fi \usepackage{ulem}\usepackage{fancybox}\usepackage{xifthen}\usepackage{forarray}\usepackage{xstring}\usepackage{stringstrings}\def\StackCreate#1#2#3{\expandafter\def\csname#1\endcsname{#3}\expandafter\def\csname#1Push\endcsname##1{\expandafter\edef\csname#1\endcsname{##1#2\csname#1\endcsname}}\expandafter\def\csname TopAux#1\endcsname ##1#2##2#3{##1}\expandafter\def\csname#1Top\endcsname{\expandafter\expandafter\expandafter\expandafter\expandafter\expandafter\csname TopAux#1\endcsname\csname#1\endcsname}\expandafter\def\csname PopAux#1\endcsname ##1#2##2#3##3#2{\expandafter\def\csname##3\endcsname{##2#3}}\expandafter\def\csname#1Pop\endcsname{\expandafter\expandafter\expandafter\expandafter\expandafter\expandafter\csname PopAux#1\endcsname\csname#1\endcsname#1#2}}\def\GetAfterColonAux#1:#2;{#2}\def\GetAfterColon#1{\IfBeginWith{#1}{:}{\GetAfterColonAux#1;}{#1}}\usepackage{aliascnt}\usepackage[shortlabels,inline]{enumitem}\setenumerate[1]{leftmargin=5.5ex}\setitemize[1]{leftmargin=5.5ex}\SetEnumitemKey{noindent}{leftmargin=0ex, itemindent=5ex, align=right, itemsep=1ex }\newcommand\NOPAGENUMBER[1]{}\usepackage{everypage}\newcommand\AddPrivateToMargin[1]{\AddEverypageHook{\tikz[overlay,remember picture]{\node at ($(current page.west)+(1.5,0)$) [rotate=90] {\textcolor{orange}{\vbox{\hrule width \the\textwidth height 0.5pt} \textcolor{black}{#1}\ \vbox{\hrule width 40em height 0.5pt}}}; }}}\newcommand\AddLongversionToMargin[1]{\AddEverypageHook{\tikz[overlay,remember picture]{\node at ($(current page.west)+(2,0)$) [rotate=90] {\textcolor{\LongColor}{\vbox{\hrule width \the\textwidth height 0.5pt} \textcolor{black}{#1}\ \vbox{\hrule width 40em height 0.5pt}}}; }}}\newcommand\AddOldversionToMargin[1]{\AddEverypageHook{\tikz[overlay,remember picture]{\node at ($(current page.west)+(2.5,0)$) [rotate=90] {\textcolor{\OldColor}{\vbox{\hrule width \the\textwidth height 0.5pt} \textcolor{black}{#1}\ \vbox{\hrule width 40em height 0.5pt}}}; }}}\newcommand\AddLineToMargin[3]{\AddEverypageHook{\tikz[overlay,remember picture]{\node at ($(current page.west)+(#2,0)$) [rotate=90] {\textcolor{#1}{\vbox{\hrule width \the\textwidth height 0.5pt} \textcolor{black}{#3}\ \vbox{\hrule width 40em height 0.5pt}}}; }}}\IfFileExists{mathabx.sty}{}{}\usepackage{amsfonts}\usepackage{amssymb}\usepackage{stmaryrd}\usepackage{amsmath}\usepackage{amsthm}\usepackage{dsfont}\IfFileExists{mbboard.sty}{\usepackage{mbboard}}{}\usepackage{mathrsfs}\usepackage{twcal}\usepackage{accents}\usepackage{bm}\usepackage[T1]{fontenc}\usepackage[latin1]{inputenc}\ifpdf \usepackage[pdftex,usenames,x11names]{xcolor}\else \usepackage[dvips,usenames,x11names]{xcolor}\fi \usepackage[pdftex]{graphicx}\usepackage[all]{xy}\ifdefined\loadTIKZ \loadTIKZ \def\TIKZlabel#1{}\else\fi \StackCreate{ColoR}{;}{?}\ColoRPush{black}\newcommand {\notion}[2][]{\def\temp{#1}\ifmmode #2 \ifx \temp\empty \index{$#2$}\else \index{$#1$}\fi \else {\bf #2}\ifx \temp\empty \index{#2}\else \index{#1}\fi \fi }\ifpdf \usepackage[pdftex,linktocpage,pagebackref,breaklinks]{hyperref}\hypersetup{colorlinks=true,allcolors=Green, linkcolor=DarkGreen, citecolor=violet, urlcolor=blue, runcolor=red, filecolor=cyan }\else \usepackage[hypertex,linktocpage,pagebackref]{hyperref}\fi \def\UndefinedRef#1{\LARGE\bfseries\color{red} ??#1??}\makeatletter \def\@setref#1#2#3{\ifx#1\relax \protect\G@refundefinedtrue \nfss@text{\reset@font\UndefinedRef{#3}}\@latex@warning{Reference `#3' on page \thepage \space undefined }\else \expandafter\Hy@setref@link#1\@empty\@empty\@nil{#2}\fi }\makeatother \usepackage{xr-hyper}\newcommand{\refX}[2]{\IfBeginWith{#1}{:}{\ref{\GetAfterColonAux#1;-#2}}{\cite[\ref{#1-#2}]{#1}}}\newcommand\pr{\begin{proof}}\def\ende{\end{proof}}\newcommand\underconstruction[1]{{\ }\\ UNDERCONSTRUCTION: #1 \ende }\newtheoremstyle{LayoutVoid}{1ex}{0ex}{\normalfont}{}{\bf}{.}{1ex}{}\newcommand\stressstatement[1]{#1}\theoremstyle{plain}\swapnumbers \newcommand\maketheorem[1]{\newtheorem{#1}[theorem]{\stressstatement{#1}} \newtheorem{#1Definition}[theorem]{\stressstatement{#1 and Definition}}  }\FunctionForEach{,}{\maketheorem}{Conclusion,Conjecture,Corollary,Fact,Facts,Lemma,Observation,Observations,Proposition,Reminder,Scholium,Summary,Theorem}\theoremstyle{definition}\theoremstyle{remark}\FunctionForEach{,}{\maketheorem}{Convention,Counterexample,Discussion,Example,Examples, Exercise,Exercises,Explanation,Notation,Project,Projects,Question,Questions,Remark,Remarks,Strategy,Warning}\theoremstyle{LayoutVoid}\numberwithin{equation}{section}\newcommand{\labelon}[1]{\marginpar{#1}}\newcommand{\labelx}[1]{{\def\temp{#1}\ifx\temp\empty\else \label{#1}\labelon{#1}\fi}}\def\GetAfterColon#1:#2;;{#2}\def\GetAfterPlus#1+#2;;{#2}\newenvironment{FACT}[2]{\IfBeginWith{#1}{:}{\def\tempFactName{void}\def\tempFreeTitle{\GetAfterColon#1;;\ }}{\IfBeginWith{#1}{+}{\def\tempFactName{voidTheorem}\def\tempFreeTitle{\GetAfterPlus#1;;\ }}{\def\tempFactName{#1}\def\tempFreeTitle{}}}\def\tempfacT{\end{\tempFactName}}\begin{\tempFactName}\labelx{#2}\textup{\textbf{\tempFreeTitle}}\capitalize[q]{#1}\caselower[q]{#1}}{\tempfacT}\catcode`\=13 \def{+}\newcommand\assigncharacter[1]{\expandafter\newcommand\csname #1\endcsname{\mathds{#1}}}\FunctionForEach{,}{\assigncharacter}{A,B,C,D,E,F,G,I,J,K,M,N,Q,R,T,U,V,W,X,Y,Z}\renewcommand\assigncharacter[1]{\expandafter\newcommand\csname C#1\endcsname{\mathcal{#1}}}\FunctionForEach{,}{\assigncharacter}{A,B,C,D,E,F,G,H,I,J,K,L,M,N,O,P,Q,R,S,T,U,V,W,X,Y,Z}\renewcommand\assigncharacter[1]{\expandafter\newcommand\csname D#1\endcsname{\mathfrak{#1}}}\FunctionForEach{,}{\assigncharacter}{a,b,c,d,e,f,g,h,i,j,k,l,m,n,o,p,q,r,s,t,u,v,w,x,y,z,A,B,C,D,E,F,G,I,K,L,M,N,O,P,Q,R,S,T,U,V,W,X,Y,Z} \renewcommand{\SS}{\mathscr{S}}\renewcommand\assigncharacter[1]{\expandafter\newcommand\csname S#1\endcsname{\mathscr{#1}}}\FunctionForEach{,}{\assigncharacter}{A,B,C,D,E,F,G,H,I,J,K,L,M,N,O,P,Q,R,T,U,V,W,X,Y,Z}\def\NewFont#1#2#3#4#5{\expandafter\font\csname #1display\endcsname =#1 at #2 \expandafter\font\csname #1normal\endcsname =#1 at #3 \expandafter\font\csname #1script\endcsname =#1 at #4 \expandafter\font\csname #1scriptscript\endcsname =#1 at #5 }\def\NewFontLetter#1#2{{\mathchoice {{\expandafter\hbox{\csname #1display\endcsname\char"#2}}}{{\expandafter\hbox{\csname #1normal\endcsname\char"#2}}}{{\expandafter\hbox{\csname #1script\endcsname\char"#2}}}{{\expandafter\hbox{\csname #1scriptscript\endcsname\char"#2}}}}}\NewFont{pxsyc}{9.00pt}{8.00pt}{7.00pt}{6.00pt}\NewFont{pxsya}{9.00pt}{8.00pt}{7.00pt}{6.00pt}\NewFont{p1xr}{10.00pt}{9.00pt}{8.00pt}{7.00pt}\NewFont{MnSymbolC10}{10.00pt}{9.00pt}{8.00pt}{7.00pt}\NewFont{MnSymbolD10}{12.00pt}{11.00pt}{10.00pt}{9.00pt}\NewFont{MnSymbolF10}{12.00pt}{11.00pt}{10.00pt}{9.00pt}\renewcommand{\bigcup}{\mathop{\NewFontLetter{MnSymbolF10}{1F}}}\NewFont{manfnt}{12.00pt}{11.00pt}{10.00pt}{9.00pt}\NewFont{favmr7y}{12.00pt}{11.00pt}{10.00pt}{9.00pt}\newcommand{\dd}{\mathrm{d}}\newcommand{\bt}{\textit{\textbf{t\,}}}\catcode95=12 \catcode95=8 \newcommand\bdl{{\ifmmode \mathrm{bdlat}\else {bounded distributive lattice}\fi}} \newcommand{\st}{{\ \vert \ }}\let\temp\phi \let\phi\varphi \let \varphi\temp \let\temp\theta \let\theta\vartheta \let \vartheta\temp \let\eps\varepsilon  \let\0\emptyset \newcommand{\into}{\hookrightarrow}\newcommand{\onto}{\twoheadrightarrow}\newcommand{\lra}{\longrightarrow}\newcommand{\Ra}{\Rightarrow}\newcommand{\La}{\Leftarrow}\newcommand{\alg}{\mathrm{alg}}\newcommand{\mal}{\cdot}\newcommand{\monthname}[1]{\ifcase#1 \or January \or February \or March \or April \or May \or June \or July \or August \or September \or October \or November \or December \fi}\newcommand\LongColor{teal}\newcommand\OldColor{gray}\newcommand\LongStart{\ColoRPush{\LongColor}\color{\ColoRTop}[BEGIN LONG VERSION]}\newcommand\LongEnd{[END LONG VERSION]\ColoRPop\color{\ColoRTop}}\newcommand\OldStart{\ColoRPush{\OldColor}\color{\ColoRTop}[BEGIN OLD VERSION]}\newcommand\OldEnd{[END OLD VERSION]\ColoRPop\color{\ColoRTop}}\newcommand\COL{\ifmmode\colon\else :\ \fi}\newcommand\kat[1]{{\tt #1}}\newcommand{\Claim}[1]{\textit{Claim #1.}\ } \renewcommand{\mod}{{\operator{\,mod\,}}}\newcommand\operator[1]{\mathop{\operatorname{#1}}\nolimits} \newcommand{\id}{\operator{id}}\newcommand{\Der}{\operator{Der}}\newcommand{\DCF}{{\rm DCF}}\newcommand{\CODF}{{\rm CODF}}\newcommand{\Gal}{\operator{Gal}}\newcommand{\qf}{\operator{qf}}\newcommand{\Hom}{\operator{Hom}}\newcommand{\Ker}{\operator {Ker}}\newcommand{\tp}{\operator{tp}}\newcommand{\ord}{\operator{ord}}\newcommand{\Spec}{\operator{Spec}}\newcommand{\Jet}{\operator{Jet}}\newcommand\DIRroot{../../../../../../../../../../../}\definecolor{Green}{rgb}{0.00,0.50,0.00}\definecolor{DarkGreen}{rgb}{0.00,0.40,0.00}\definecolor{grey}{rgb}{0.40,0.40,0.40} \renewcommand\textcolor[2]{\ColoRPush{#1}\color{\ColoRTop}#2\ColoRPop\color{\ColoRTop}}\IfFileExists{C:/wb/System64/WinBatch.exe}{}{}\let\temp\theta \let\theta\vartheta \let \vartheta\temp \def\uc{{\rm UC} }\def\Ind{\setbox0=\hbox{$x$}\kern\wd0\hbox to 0pt{\hss$\mid$\hss} \lower.9\ht0\hbox to 0pt{\hss$\smile$\hss}\kern\wd0}\def\Notind{\setbox0=\hbox{$x$}\kern\wd0\hbox to 0pt{\mathchardef \nn=12854\hss$\nn$\kern1.4\wd0\hss}\hbox to 0pt{\hss$\mid$\hss}\lower.9\ht0 \hbox to 0pt{\hss$\smile$\hss}\kern\wd0}\def \DCF {\operatorname{DCF}}\def \UC {\operatorname{UC}}\def \jet {\operatorname{Jet}}\def \pol {\operatorname{pol}}\def \SS {\mathcal S}\def \V{\mathcal V}\def \DD {\mathcal D}\newcommand \ec {e.c.}\newcommand\ev{\mathrm{e\hspace{-.5pt}v}}\renewcommand{\labelon}[1]{}\renewcommand\LongStart{\ColoRPush{\LongColor}\color{\ColoRTop}}\renewcommand\LongEnd{\ColoRPop\color{\ColoRTop}}\ifprivate \AddPrivateToMargin{private version} \fi \iflongversion \AddLongversionToMargin{long version} \fi \ifoldversion \AddOldversionToMargin{old version included} \fi 
\begin{document} \title{Differentially Large Fields} \author{Omar Le\'on S\'anchez} \address{Omar Le\'on S\'anchez, The University of Manchester\\ Department of Mathematics\\ Oxford Road \\ Manchester, M13 9PL, UK} \email{omar.sanchez@manchester.ac.uk} \author{Marcus Tressl} \address{Marcus Tressl, The University of Manchester\\ Department of Mathematics\\ Oxford Road \\ Manchester, M13 9PL, UK \newline Homepage: \homepageTressl} \email{\emailTressl} \date{\today}  \begin{abstract} We introduce the notion of \textit{differential largeness} for fields equipped with several commuting derivations (as an analogue to largeness of fields). We lay out the foundations of this new class of "tame" differential fields. We state several characterizations and exhibit plenty of examples and applications. Our results strongly indicate that differentially large fields will play a key role in differential field arithmetic. For instance, we characterise differential largeness in terms of being existentially closed in their power series field (furnished with natural derivations), we give explicit constructions of differentially large fields in terms of iterated powers series, we prove that the class of differentially large fields is elementary, and we show that differential largeness is preserved under algebraic extensions, therefore showing that their algebraic closure is differentially closed. \end{abstract} \maketitle \setcounter{tocdepth}{1} \tableofcontents \section{Introduction} \noindent Recall that a field $K$ is called \textit{large} (or \textit{ample}) if every irreducible variety defined over $K$ with a smooth $K$-rational point has a Zariski-dense set of $K$-rational points. Equivalently, every variety defined over $K$ that has a $K((t))$-rational point, also has a $K$-rational point. Large fields constitute one of the widest classes of \textit{tame fields}; namely, every class of fields that serves as a \textit{locality}, in the sense that universal local-global principles hold, consists entirely of large fields, cf. \cite{BaSFeh2013,Pop2013}. For example, all local fields are large and so are pseudo-classically closed fields (like PAC or PRC fields), the field of totally real numbers, as well as the quotient field of any local henselian domain \cite{Pop2010}. On the other hand, number fields and algebraic function fields are not large by Faltings' theorem and its function field version. \par One of the most remarkable Galois-theoretic applications of large fields, due to Pop \cite{Pop1996}, states every finite split embedding problem over large fields has proper regular solutions. In particular, the regular inverse Galois problem is solvable over all large fields. Pop's work (and the work of many others) demonstrates that ``over large fields one can do a lot of interesting mathematics". For instance, large fields have been widely used to tackle long-standing problems in field arithmetic: inverse Galois theory, torsors of finite groups, elementary theory of function fields, extremal valued fields, to name a few. We refer the reader to Pop's survey \cite{Pop2013} for earlier and current developments on the subject, and to \cite{BaSFeh2013} for a list of open problems. \par In this paper we introduce the notion of \textit{differential largeness} in the class of differential fields of characteristic 0 in several commuting derivations. We lay out the foundations of this new and exciting class of "tame" differential fields, prove several characterisations, and exhibit plenty of examples and applications. In order to give the definition of a differentially large field we need one piece of terminology. We say that a field $K$ is \textit{existentially closed} (e.c.) in $L$ if every variety defined over $K$ that has an $L$-rational point, also has a $K$-rational point. Hence a field is large just if it is \ec\ in its Laurent series field. Similarly, a differential field $K$ (of characteristic 0 throughout, in $m\geq 1$ commuting derivations) is \ec\ in a differential field extension $L$ if every differential variety defined over $K$ that has an $L$-rational differential point, also has a $K$-rational differential point. (See \ref{ECbasic} for other characterizations of this property.) \par A differential field is \textit{differentially large} if it is large as a pure field and for every differential field extension $L/K$, if $K$ is \ec\ in $L$ as a field, then it is \ec\ in $L$ as a differential field. For example differentially closed fields (aka \textit{constrainedly closed} in Kolchin terminology) and closed ordered differential fields in the sense of \cite{Singer1978a} are differentially large. \par In Theorem \ref{CharDlargeI}, we establish several equivalent formulations of differential largeness that justify why indeed this is the right differential analogue of largeness. For instance, we characterise them in terms of differential varieties having a Kolchin-dense set of rational points as long as they have suitable ``smooth" rational points. In addition, we prove (in analogy to the characterisation of largeness in terms of being e.c. in its Laurent series field) that a differential field $K$ is differentially large just if it is e.c. in its power series field $K((t_1,\dots,t_m))$ as differential fields. The derivations on the power series field are given by the unique commuting derivations $\delta_1,\dots,\delta_m$ extending the ones on $K$ that are compatible with infinite sums and satisfy $\delta_i(t_j)=\frac{d t_j}{dt_i}$. \par A key tool in establishing our formulations of differential largeness (and further results) is the introduction of a twisted version of the classical Taylor morphism associated to a ring homomorphism $\phi:A\lra B$ for a given differential ring $A$. We explain this briefly in the case of one derivation $\delta$. Recall that the Taylor morphism $T_\phi(a)=\sum_{k\geq 0}\frac{\phi(\delta^k(a))}{k!}t^k$ defines a differential ring homomorphism $(A,\delta)\lra (B[[t]],\frac{\dd}{\dd t})$. Typically this is applied when $A$ is a differential $K$-algebra for a differential field $K$ and $\phi $ is a (not necessarily differential) $K$-algebra homomorphism $A\lra K$ (so $B=K$). If the derivation on $K$ is trivial, then $T_\phi$ is in fact a (differential) $K$-algebra homomorphism and in this context it was used by Seidenberg for example, to establish his embedding theorem for differential fields into meromorphic functions. However, if the derivation on $K$ is not trivial, then $T_\phi$ is not a $K$-algebra homomorphism, i.e, it is not an extension of $\phi$. On the other hand, $T_\phi$ can be ``twisted" in order to obtain a natural differential $K$-algebra homomorphism $T_\phi^*:(A,\delta)\lra (K[[t]],\partial)$, where $\partial $ is the natural derivation extending the given one on $K$ and satisfying $\partial (t)=1$. This is established in \ref{twistedTaylor}, where we use it to derive the following result that seems to be of independent interest; for instance, in the analysis of formal solutions to PDEs, see \cite{Sanz2017}, and is deployed in most parts of this article (in the more general form \ref{DiffPointFromPointInECext}. \par \medskip \par \noindent {\bf Theorem.}\textit{ Let $(K, \delta)$ be a differential field of characteristic zero that is large as a field and let $(S,\delta)$ be a differentially finitely generated $K$-algebra. If there is a $K$-algebra homomorphism $S\to L$ for some field extension $L/K$ in which $K$ is \ec\ (as a field), then there is a differential $K$-algebra homomorphism $(S,\delta)\to (K[[\bt]],\partial)$. } \par \medskip\noindent Differentially large fields will play a very similar role in differential field arithmetic to that played by large fields in field arithmetic (of characteristic 0). The principal indicators for this are established in this paper (in Sections \ref{sectionUCalgebraically} and \ref{further}). We show that \begin{enumerate}[(a)] \item A differential field $K$ is differentially large if and only if it is existentially closed in its power series field $K((t_1,\ldots,t_m))$ furnished with $m$ natural derivations extending those on $K$ satisfying $\partial_i(t_j)=\frac{dt_j}{dt_i}$. See \ref{CharDlargeI}. \item Every large field equipped with commuting derivations has an extension to a differentially large field $L$ such that $K$ is \ec\ in $L$ on the level of pure fields. See \ref{UCrecall}. \item Differentially large fields are first order axiomatisable (see \ref{DLisUC} and also \ref{DiffLargeGeomAx} for a concrete algebro-geometric description), and the elimination theory of the underlying field transfers to the differential field; see \ref{UCrecall}. \item Differential largeness is preserved under algebraic extensions. Thus, the algebraic closure of a differentially large field is differentially closed. This provides many new differential fields with minimal differential closures. see \ref{algextdifflarge}. \item Differentially large fields (and differentially closed fields) can be produced by iterated power series constructions. See \ref{ConstructionPowerSeries}. \item The existential theory of the class of differentially large fields is the existential theory of the differential field $\mathbb Q((\bt_1))((\bt_2))$ equipped with its natural derivations, see \ref{existentialresult}. \item Differentially large fields are Picard-Vessiot closed, see \ref{DiffLargeImpliesPVclosed}. \item Connected differential algebraic groups defined over differentially large fields have a Kolchin-dense set of rational points, \ref{KolDenseInDiffGroups}. \item A differentially large field is PAC (at the field level) if and only it is pseudo differentially closed, see \ref{psfDCF}. \end{enumerate} \par \noindent Large fields have also made an appearance in the (inverse) Picard-Vessiot theory of linear ordinary differential equations. In \cite{BHHP2020}, it is shown that if $K$ is a large field of infinite transcendence degree, then every linear algebraic group over $K$ is a Picard-Vessiot group over $(K(x), \frac{\dd}{\dd x})$. We envisage that differentially large fields will make a similar appearance in the Parameterised Picard-Vessiot theory and its differential (constrained) coholomogy. The first application in this direction already appears in a paper of the first author with A. Pillay \cite{LeonPillay2019} using an earlier draft of the present paper. They show that if an ordinary differential field $(K,\delta)$ is differentially large and bounded as a field (i.e., has finitely many extensions of degree $n$, for each $n\in\N)$, then for any linear differential algebraic group $G$ over $K$ the differential Galois cohomology $H_\delta^1(K,G)$ is finite. This can be thought of as a differential analogue of the classical result of Serre stating that if a field $K$ is bounded then the Galois cohomology $H^1(K,G)$ is finite for any linear algebraic group over $K$. \section{Preliminaries} \noindent All rings and algebras in this article are assumed to be commutative and unital. We also assume that all our fields are of characteristic zero. \par In this section we briefly summarize the key notions and terminology, mostly from differential algebra, that we will freely use throughout the paper (especially in Section \ref{sectionUCalgebraically} where we give several equivalent formulations of differential largeness). We make a few remarks on the notion of existentially closed differential ring extensions, we recall the Structure Theorem for finitely generated differential algebras, and give a quick review of jets and prolongation spaces. \par \medskip\noindent Recall that a derivation on a ring $R$ is an additive map $\delta:R\to R$ satisfying the Leibniz rule $$\delta(rs)=\delta(r) s+r\delta(s) \quad \text{ for all } r,s\in R.$$ Throughout a differential ring $R=(R,\Delta)$ is a ring $R$ equipped with a distinguished set of \textit{commuting} derivations $\Delta=\{\delta_1,\dots,\delta_m\}$. We also allow the case when $m=0$, in which case we are simply talking of rings with no additional structure. \par Given a differential ring $R$, a differential $R$-algebra $A$ is an $R$-algebra equipped with derivations $\Delta=\{\delta_1,\dots,\delta_m\}$ such that the structure map $R\to A$ is a differential ring homomorphism. If $L$ is another differential $R$-algebra which is also a field, then an \notion{$L$-rational point} of $A$ is a differential $R$-algebra homomorphism $A\lra L$. This terminology is in line with the standard language of algebraic geometry, where $A$ is thought of as $R\{x_1,\ldots,x_n\}/I$, with $I$ a differential ideal of the differential polynomial ring $R\{x_1,\ldots,x_n\}$, and the differential $R$-algebra homomorphisms $A\lra L$ are coordinate free descriptions of the common differential zeroes $a\in L^n$ of the polynomials from $I$ (via evaluation at $a$). \par For the basics in differential algebra, such as differential field extensions and differentially closed fields (also called constrainedly closed which is the differential analogue of algebraically closed), we refer the reader to the excellent book of Kolchin \cite{Kolchi1973}.\looseness=-1 \par \par \bigskip \begin{FACT}{:Existentially closed extensions.}{ECbasic} Fix $m\geq 0$. Let $B=(B,\delta_1,\ldots,\delta_m)$ be a differential ring and let $A$ be a differential subring of $B$. (If $m=0$, $B$ is just a ring and $A$ is a subring.) Then $A$ is said to be \notion{existentially closed (e.c.) in $B$} if for every $n\in\N$ and all finite collections $\Sigma,\Gamma\subseteq A\{x_1,\ldots,x_n\}$ of differential polynomials in $m$ derivations and $n$ differential variables, if there is a common solution in $B^n$ of $P=0\ \&\ Q\neq 0$ ($P\in\Sigma,Q\in\Gamma$), then such a solution may also be found in $A^n$. \par We are mainly interested in the case when $A=K$ is a differential field and in this case we will use the following properties (in the case $m=0$, differentially finitely generated, differential field, etc. should be understood as finitely generated, field, etc., and differentially closed field should be understood as algebraically closed field and Kolchin topology should be understood as Zariski topology). We make heavy use of the following properties: \begin{enumerate}[(i),itemsep=1ex] \item\labelx{ECbasicInDomainToField} If $K$ is \ec\ in $B$, then one easily checks that $B$ is a domain and that $K$ is also \ec\ in $\qf(B)$. \item\labelx{ECbasicFieldField} If $B$ is also a differential field, then $K$ is \ec\ in $B$ if and only if every differentially finitely generated $K$-algebra $S$ that possesses a differential point $S\lra B$, also possesses a differential point $S\lra K$. The reason is that when $B$ is a field then the inequalities $Q\neq 0$ in the definition of existentially closed above may be replaced by the equality $y\mal Q(x)-1=0$, where $y$ is a new variable. \item\labelx{ECbasicECasKdense} If $B$ is a differentially finitely generated $K$-algebra then the following are equivalent. \begin{enumerate}[(a)] \item $K$ is \ec\ in $B$. \item $B$ is a domain and for each $b\in B$, if $f(b)=0$ for every differential $K$-rational point $f:B\lra K$, then $b=0$.\footnote{In other words, in the subspace of $\Spec(B)$ consisting of differential prime ideals, the set of maximal and differential ideals with residue field $K$, is dense.} (In particular $B$ has a differential $K$-rational point.) We refer to this property as \notion{$B$ has a Kolchin dense set of differential $K$-rational points}. \item For all $n\in\N$, each differential prime ideal $\Dp$ of $K\{x\}$, $x=(x_1,\ldots,x_n)$ with $B\cong_KK\{x\}/\Dp$ and each differential field $L$ containing $K$, the set $V_K=\{a\in K^n\st \Dp(a)=0\}$ is dense in $V_L=\{a\in L^n\st \Dp(a)=0\}$ for the \notion{Kolchin topology} of $L^n$ (having zero sets of differential polynomials from $L\{x\}$ as a basis of closed sets). \item There is some $n\in\N$, a differential prime ideal $\Dp$ of $K\{x\}$, $x=(x_1,\ldots,x_n)$ with $B\cong_KK\{x\}/\Dp$ and a differentially closed field $M$ containing $K$ such that the set $V_K$ is dense in $V_M$ for the $K$-Kolchin topology of $M^n$ (having the zero sets in $M^n$ of differential polynomials from $K\{x\}$ as a basis of closed sets). \end{enumerate} \par \noindent \textit{Proof of} (iii). We may assume that $B$ is a domain throughout and write $B=K\{x\}/\Dp$, $x=(x_1,\ldots,x_n)$. The arguments below go through for any choice of these data. By the differential basis theorem there is some finite $\Phi\subseteq K\{x\}$ such that $\Dp$ is the radical differential ideal $\sqrt[d]\Phi$ generated by $\Phi$. For a differential field $L$ containing $K$ we write $V_L=\{a\in L^n\st \Dp(a)=0\}$ and $I_L=\{Q\in L\{x\}\st Q|_{V_L}=0\}$. If $L$ is differentially closed, then the differential Nullstellensatz \cite[Ch. IV, section 3, Theorem 2, p. 147]{Kolchi1973} says $I_L=\sqrt[d]\Dp$ (in $L\{x\}$).\looseness=-1 \par \medskip\noindent (a)$\Ra $(b). If $K$ is \ec\ in $B$ and $b\in B\setminus \{0\}$, then take $Q\in K\{x\}$ with $b=Q(x+\Dp)$. Since in $B$ we have a solution of $\Phi=0\ \&\ Q\neq 0$, there is also a solution $a\in K^n$ and evaluation $K\{x\}\lra K$ at $a$ factors through a differential $K$-rational point $B\lra K$ that is non-zero at $b$. \par \medskip\noindent (b)$\Ra $(c). Let $M$ be a differentially closed field containing $L$ such that the fixed field of the group of differential $K$-automorphisms of $M$ is $K$ (for example $M$ could be a sufficiently saturated differentially closed field or, in Kolchin's terminology, a universal differential extension of $K$). It suffices to show that $V_K$ is dense in $V_M$ for the Kolchin topology of $M^n$. Let $W$ be the closure of $V_K$ in $M^n$ for the Kolchin topology of $M^n$ and let $J=\{Q\in M\{x\}\st Q|_W=0\}$. Then every differential $K$-automorphism of $M$ fixes $W$ set wise and so also fixes $J$ set wise. Using our assumption on $M$ we see that the differential field of definition of $J$ is contained in $K$; hence the differential ideal $J$ is generated as an ideal by $J\cap K\{x\}$. On the other hand we have $I_M\cap K\{x\}=\Dp$. As $W\subseteq V_M$ we get $\Dp \subseteq I_M\subseteq J$ and we claim that $\Dp=J\cap K\{x\}$. Take $P\in J\cap K\{x\}$. Then $P$ vanishes on $V_K$, which says that the element $P+\Dp\in B$ is mapped to $0$ by all differential $K$-rational points of $B$. By (b), this implies $P+\Dp=0$ in $B$, in other words  $P\in \Dp$. We have shown that $\Dp=I_M\cap K\{x\}=J\cap K\{x\}$, which implies $W=V_M$ as required. \par \smallskip\noindent (c)$\Ra $(d) is trivial. \par \smallskip\noindent (d)$\Ra $(a). Let $\Sigma=\{P_1,\ldots,P_s\},\Gamma\subseteq K\{y\}$ be finite, $y=(y_1,\ldots,y_r)$ and assume there is some $c\in B^r$ with $\Sigma(c)=0\ \&\ \Gamma(c)\neq 0$. We need to find some $a\in K^r$ with $\Sigma(a)=0\ \&\ \Gamma(a)\neq 0$. Since $B$ is a domain we may assume that $\Gamma=\{Q(y)\}$ is a singleton. We write $c_i=H_i(x+\Dp)$ with $H_i\in K\{x\}$ and $H=(H_1,\ldots,H_r)$. Then $\Sigma(c)=0\ \&\ Q(c)\neq 0$ means $P_1(H),\ldots,P_s(H)\in \Dp$ and $Q(H)\notin \Dp$. Since $M$ is differentially closed and $Q(H)\notin \Dp=K\{x\}=I_M\cap K\{x\}$ there is some $d\in V_M$ with $Q(H(d))\neq 0$. By (d) and because $Q(H)\in K\{x\}$, there is some $b\in V_K$ with $Q(H(b))\neq 0$. Since $P_i(H)\in \Dp$ we also know $P_i(H(b))=0$. Hence the tuple $a=(H_1(b),\ldots,H_r(b))\in K^r$ solves the given system. \hfill$\diamond$ \item\labelx{ECbasicECsmooth} If $m=0$ and $B$ is a finitely generated $K$-algebra then $K$ is \ec\ in $B$ if and only if $B$ is a domain and the set of smooth $K$-rational points of $B$ is Zariski dense in the $L$-rational points for any field $L$ containing $K$. This is a statement in classical algebraic geometry (using the formulation (c) of \ec\ in \ref{ECbasicECasKdense}). If in addition $K$ is a large field, then $K$ is \ec\ in $B$ if and only if $B$ is a domain that has a smooth $K$-rational point. \par \end{enumerate} \par \noindent If $B$ is a differential $K$-algebra we will say that \notion{$K$ is existentially closed in $B$ as a field} if it is \ec\ in $B$ when we forget about the derivations; hence if the above condition holds true for systems $\Sigma,\Gamma$ of ordinary (non-differential) polynomials. If we want to emphasize that the derivations are to be taken into account we say \notion{$K$ is existentially closed in $B$ as a differential field}. \end{FACT} \bigskip \begin{FACT}{:Structure theorem for finitely generated differential algebras.}{structuretheorem} Let $K$ be a differential field (of characteristic 0) and let $S$ be a differential $K$-algebra that is differentially finitely generated and a domain. Then, by \cite{Tressl2002}, there are $K$-subalgebras $A,P$ of $S$ and an element $h\in A\setminus \{0\}$ such that $A$ is a finitely generated $K$-algebra, $P$ is a polynomial $K$-algebra and the natural homomorphism $A_h\otimes_KP\lra S_h$ given by multiplication, is an isomorphism. (Note that in general neither $A_h$ nor $P$ is differential.) \end{FACT} \begin{FACT}{Definition}{defnComposite} Let $K$ be a differential field. We call a differential $K$-algebra $S$ \notion[]{composite} if $S$ is finitely generated as a differential $K$-algebra and a domain, which possess (not necessarily differential) $K$-subalgebras $A,P$ such that \begin{enumerate}[(a)] \item $A$ is a finitely generated $K$-algebra and $P$ is a polynomial $K$-algebra, and \item the natural map $A\otimes_KP\lra S$ is an isomorphism. \end{enumerate} \end{FACT} \begin{FACT}{Corollary}{StructureThmConsequence} Let $K$ be a differential field and let $S$ be a finitely generated differential $K$-algebra. Let $f:S\lra L$ be a differential $K$-algebra homomorphism to some differential field extension $L$ of $K$. \begin{enumerate}[(i)] \item There is a differential $K$-subalgebra $S_0\subseteq L$ that is composite and contains the image of $f$. \item If $K$ is \ec\ in $L$ as a field, then there is a $K$-algebra homomorphism $S\lra K$. \end{enumerate} \end{FACT} \begin{proof} (i). Let $\Dp$ be the kernel of $f$. Then $S/\Dp$ is again a differentially finitely generated $K$-algebra and so we may assume that $\Dp=0$ and $S\subseteq L$. By the Structure Theorem \ref{structuretheorem}, $S_h\cong_K A_h\otimes_K P$ where $A$ is a finitely generated $K$-subalgebra of $S$, $h\in A$, and $P$ is a polynomial $K$-algebra, $P\subseteq S$. As $S\subseteq L$, we have $A_h\subseteq L$. Hence we may take $S_0=S_h$. \par \smallskip\noindent (ii). Take $S_0$ as in (i) and $A,P$ for $S_0$ as in \ref{defnComposite}. Since $A$ is a finitely generated $K$-subalgebra of $L$ and $K$ is \ec\ in $L$ as a field, there is a $K$-algebra homomorphism $A\to K$. Since $P$ is a polynomial $K$-algebra there is also a $K$-algebra homomorphism $P\lra K$. Hence by the universal property of the tensor product there is a $K$-algebra homomorphism $S\lra K$. \end{proof} \par \bigskip\noindent \begin{FACT}{:Differential Varieties, Jets and Prolongations.}{diffalggeometry}We recall the basic objects of differential algebraic geometry in the sense of Kolchin \cite{Kolchi1973}, and the constructions of jets and prolongations. Some parts are notationally heavy but we try to only introduce those that we will need (and freely use) in coming sections. \end{FACT} We work inside a (sufficiently saturated or universal) differentially closed field $(\U,\Delta)$, and $K$ denotes a differential subfield of $\U$. A \notion{Kolchin-closed} subset of $\U^n$ is the common zero set of a set of differential polynomials over $\U$ in $n$ differential variables; such sets are also called \notion{affine differential varieties}. If the definining polynomials can be chosen with coefficients in $K$ we will say the set is \notion[]{defined over $K$}. \par By a \notion{differential variety} $V$ we mean a topological space which has as finite open cover $V_1,\dots,V_s$ with each $V_i$ homeomorphic to an affine differential variety (inside some power of $\U$) such that the transition maps are regular as differential morphisms; see \cite[Chap. 1, section 7]{LeoSan2013}. We will say that the differential variety is over $K$ when all objects and morphisms can be defined over $K$. This definition also applies to our use of algebraic varieties, replacing Kolchin-closed with Zariski-closed in powers of $\U$ (recall that $\U$ is algebraically closed and a universal domain for algebraic geometry in Weil's ``foundations'' sense). \begin{FACT}{void}{} We fix integers $n>0$ and $r\geq 0$, and set $$ \Gamma_n(r) = \{(\xi,i) \in \mathbb N^m\times\{1,\dots,n\} \st \sum_{i=1}^m \xi_i \leq r\}. $$ \par \noindent The \notion[]{$r$-th nabla map} $\nabla_r:\U^n\to \U^{\alpha(n,r)}$ with $\alpha(n,r):=|\Gamma_n(r)|=n\cdot\binom{r+m}{m}$ is defined by $$\nabla_r(x)= (\delta^\xi x_i:\,(\xi,i)\in \Gamma_n(r)),$$ where $x=(x_1,\dots,x_n)$ and $\delta^\xi=\delta_1^{\xi_1}\cdots\delta_m^{\xi_m}$. We order the elements of the tuple $(\delta^\xi x_i:\,(\xi,i)\in \Gamma_n(r))$ according to the canonical orderly ranking of the indeterminates $\delta^\xi x_i$; that is, \begin{equation}\label{ordercanonical} \delta^{\xi}x_i< \delta^{\zeta}x_j \iff \left(\sum \xi_k,i,\xi_1,\dots,\xi_m\right)<_{\text{lex}}\left(\sum \zeta_k,j,\zeta_1,\dots,\zeta_m\right) \end{equation} \par Let $\U_r:=\U[\epsilon_1,\dots,\epsilon_m]/(\epsilon_1,\dots,\epsilon_m)^{r+1}$ where the $\epsilon_i$'s are indeterminates, and let $e:\U\to \U_r$ denote the ring homomorphism $$x\mapsto \sum_{\xi\in\Gamma_1(r)}\frac{1}{\xi_1!\cdots\xi_m!}\; \delta^\xi(x)\; \epsilon_1^{\xi_1}\cdots\epsilon_m^{\xi_m}.$$ We call $e$ the exponential $\U$-algebra structure of $\U_r$. To distinguish between the standard and the exponential algebra structure on $\U_r$, we denote the latter by $\U_r^e$. \end{FACT} \begin{FACT}{Definition}{} Given an algebraic variety $X$ the $r$-th \notion{prolongation} $\tau X$ is the algebraic variety given by the taking the $\U$-rational points of the classical Weil descent (or Weil restriction) of $X\times_\U \U_r^e$ from $\U_r$ to $\U$. Note that the base change $V\times_\U \U_r^e$ is with respect to the exponential structure while the Weil descent is with respect to the standard $\U$-algebra structure. \end{FACT} \smallskip\noindent For details and properties of prolongation spaces we refer to \cite[\S 2]{MoPiSc2008}; for a more general presentation see \cite{MooSca2010}. In particular, it is pointed out there that the prolongation $\tau_r X$ always exist when $X$ is quasi-projective (an assumption that we will adhere to later on). A characterising feature of the prolongation is that for each point $a\in X=X(\U)$ we have $\nabla_r(a)\in \tau_r X$. Thus, the map $\nabla_{r}:X\to \tau_r X$ is a differential regular section of $\pi_r:\tau_r X\to X$ the canonical projection induced from the residue map $\U_r\to \U$. We note that if $X$ is defined over the differential field $K$ then $\tau_r X$ is defined over $K$ as well. \par In fact, $\tau_r$ as defined above is a functor from the category of algebraic varieties over $K$ to itself, and the maps $\pi_r:\tau_r X\to X$ and $\nabla_r:X\to \tau_r X$ are natural. The latter means that for any morphism of algebraic varieties $f:X\to Y$ we get \begin{equation}\label{natural} f\circ\pi_{r,X}=\pi_{r,Y}\circ\tau_r f \quad \text{ and }\quad \tau_r f\circ\nabla_{r,X}=\nabla_{r,Y}\circ f. \end{equation} If $G$ is an algebraic group, then $\tau_r G$ also has the structure of an algebraic group. Indeed, since $\tau_r$ commutes with products, the group structure is given by $$\tau_r(*):\tau_r G\times\tau_r G\to \tau_r G$$ where $*$ denotes multiplication in $G$. Moreover, by the right-most equality in \eqref{natural}, the map $\nabla_r:G\to\tau_r G$ is an injective group homomorphism. Hence, $\nabla_r(G)$ is a differential algebraic subgroup of $\tau_r G$. We will use this in \ref{KolDenseInDiffGroups} below. \par Assume that $V$ is a differential variety which is given as a differential subvariety of a quasi-projective algebraic variety $X$. We define the $r$-th jet of $V$ to be the Zariski-closure of the image of $V$ under the $r$-th nabla map $\nabla_r:X\to \tau_r X$; that is, $$\jet_r V=\overline{\nabla_r(V)}^{\operatorname{Zar}}\subseteq \tau_r X.$$ The jet sequence of $V$ is defined as $(\jet_r V:r\geq 0)$. Note that this sequence determines $V$, indeed $$V=\{a\in X: \nabla_r(a)\in \jet_r V \text{ for all $r\geq 0$}\}.$$ \par \medskip \begin{FACT}{:General Assumption.}{} Throughout we assume, whenever necessary for the existence of jets, that our differential varieties are given as differential subvarieties of quasi-projective algebraic varieties. Of course, in the affine case this is always the case. It is worth noting, as it will be used in \ref{KolDenseInDiffGroups}, that for connected differential algebraic groups this is also true. Indeed, by \cite[Corollary 4.2(ii)]{Pillay1997} every such group embeds into a connected algebraic group and the latter is quasi-projective by Chevalley's theorem. \end{FACT} \goodbreak \section{The Taylor Morphism}\label{taylor} \numberwithin{theorem}{section} \par \noindent In parallel to the characterization of large fields in terms of being e.c. in Laurent series, we will prove in \ref{CharDlargeI} that differential largeness can be characterized similarly. For this, we will make use of a \textit{twisted} Taylor morphism. In this section, we give a description of this morphism and use it to construct solutions in power series to systems of differential equations (cf. Corollary~\ref{DiffPointFromPointInECext}) \begin{FACT}{void}{SetupTaylor} Let $(A,\Delta)$ be a differential ring with commuting derivations $\Delta=\{\delta_1,\dots,\delta_m\}$. Recall that given a ring homomorphism $\phi:A\to B$ (where $B$ is a $\Q$-algebra), the Taylor morphism $T_\Delta^\phi:A\to B[[\bt]]$, where $\bt$ $=(t_1,\dots,t_m)$, is defined as \[ a\mapsto \sum_{\alpha} \frac{\phi(\delta^\alpha a)}{\alpha!}\; \bt^\alpha \] where we make use of multi-index notation. Namely, $\alpha=(\alpha_1,\dots,\alpha_m)\in\N^m$, $\alpha!=\alpha_1!\cdots\alpha_m!$, $\delta^\alpha=\delta_1^{\alpha_1}\cdots\delta_m^{\alpha_m}$, and ${\bt}^{\alpha}=t_1^{\alpha}\cdots t_m^{\alpha_m}$. It is a straightforward computation to check that $T_\Delta^\phi$ is a differential ring homomorphism \[ (A,\Delta) \to \left(B[[{\bt}]], \frac{\dd}{\dd t_1}, \dots,\frac{\dd}{\dd t_m}\right). \] \par \smallskip\noindent For every such family of commuting derivations $\Delta$ on $A$, there is a unique extension to $A[[\bt]]$ such that the derivations commute with meaningful sums and map all $t_i$'s to $0$ (cf. \cite[Chap. 0, section 13]{Kolchi1973}). We continue to denote these derivations on $A[[\bt]]$ by $\Delta=\{\delta_1,\dots,\delta_m\}$; note that they still commute with each other. We work with the derivations $\delta_i+\frac{\dd}{\dd t_i}$, for $i=1,\dots,m$, on $A[[\bt]]$; again these commute with each other. Assuming that $A$ is a $\Q$-algebra, we now study the algebraic properties of the Taylor morphism associated to the evaluation map $$\ev:A[[{\bt}]]\to A[[{\bt}]], \quad f\mapsto f(0,\dots,0).$$ For instance, we show that the map from $Der(A)$ to ring endomorphisms of $A[[t]]$ given by $\delta\mapsto T^{\ev}_{\delta+\frac{\dd}{\dd t}}$ is a monoid homomorphism when restricted to any submonoid of commuting derivations. Here the monoid structure on $Der(A)$ is just addition of derivations (and so is indeed a group), while the monoid structure on ring endomorphisms is composition. Note that as a consequence $T^{\ev}_{\delta+\frac{\dd}{\dd t}}$ is a differential ring isomorphism, because $T_{\delta+\frac{\dd}{\dd t}}^{\ev}$ has compositional inverse $T_{- \delta+\frac{\dd}{\dd t}}^{\ev}$ and $T_{\frac{\dd}{\dd t}}^{\ev}$ is the identity map on $A[[t]]$. We state all this more generally below. \end{FACT} \smallskip We first introduce some convenient notation and terminology. Let $\Delta=\{\delta_1,\dots,\delta_m\}$ and $\Omega=\{\partial_1,\dots,\partial_m\}$ be families of commuting derivations on $A$. We say that these families commute if $\delta_i$ commutes with $\partial_j$ for all $1\leq i,j\leq m$; when this is the case, we denote by $\Delta+\Omega$ the family of commuting derivations on $A$ given by $\{\delta_1+\partial_1,\dots,\delta_m+\partial_m\}$. Note that the natural extensions of $\Delta$ and $\Omega$ to $A[[\bt]]$, as discussed above, commute with the family $$\frac{\dd}{\dd \bt}:=\left\{\frac{\dd}{\dd t_1},\dots,\frac{\dd}{\dd t_m}\right\}.$$ Therefore, the family of derivations $\Delta+\Omega+\frac{\dd}{\dd \bt}$ on $A[[\bt]]$ is a commuting family. \begin{FACT}{Theorem}{propertyT} Let $A$ be a $\Q$-algebra, and let $\Delta$ and $\Omega$ be families of $m$-many commuting derivations on $A$. If $\Delta$ and $\Omega$ commute, then \begin{equation}\label{equaltaylor} T_{\Delta+\Omega+\frac{\dd}{\dd \bt}}^{\ev}\;=\; T_{\Delta+\frac{\dd}{\dd \bt}}^{\ev}\; \circ\; T_{\Omega+\frac{\dd}{\dd \bt}}^{\ev}. \end{equation} \end{FACT} \begin{proof} For $\alpha\in \N^m$ we use the multi-index notation \begin{align*} (\delta+\partial)^\alpha&=(\delta_1+\partial_1)^{\alpha_1}\cdots(\delta_m+\partial_m)^{\alpha_m},\text{ and}\\ (\delta+\partial+\frac{\dd}{\dd \bt})^\alpha&=(\delta_1+\partial_1+\frac{\dd}{\dd t_1})^{\alpha_1}\cdots(\delta_m+\partial_m+\frac{\dd}{\dd t_m})^{\alpha_m}. \end{align*} We use the product order $\leq$ on $\N^m$ given by $\beta\leq \alpha$ if and only if $\beta_i\leq \alpha_i$ for $1\leq i\leq m$). As the derivations commute, we have the usual binomial identities \begin{align*} (\delta+\partial)^\alpha&= \sum_{\beta\leq \alpha}\binom{\alpha}{\beta}\delta^\beta\partial^{\alpha-\beta} =\sum_{\beta+\gamma= \alpha}\binom{\alpha}{\beta}\delta^\beta\partial^\gamma,\text{ and}\cr (\delta+\partial+\frac{\dd}{\dd \bt})^\alpha&=\sum_{\xi\leq\alpha}\sum_{\beta+\gamma=\xi}\binom{\alpha}{\xi}\binom{\xi}{\beta}\delta^{\beta}\partial^\gamma \frac{\dd^{\alpha-\xi}}{\dd \bt} \\ &= \sum_{\beta+\gamma\leq \alpha}\binom{\alpha}{\beta+\gamma}\binom{\beta+\gamma}{\beta}\delta^\beta\partial^\gamma\frac{\dd^{\alpha-\beta-\gamma}}{\dd \bt}. \end{align*} \par \noindent Now take $f=\sum_{\xi}a_\xi \bt^{\xi}\in A[[\bt]]$. We show that both sides of equation \eqref{equaltaylor} applied to $f$ are equal to \begin{equation}\label{bothsides} \sum_{\alpha}\left(\sum_{\beta+\gamma\leq \alpha}\frac{1}{\beta!\cdot \gamma!}\; \delta^\beta\partial^\gamma(a_{\alpha-\beta-\gamma})\right)\bt^\alpha. \end{equation} \noindent We begin with the left-hand-side. By definition, the coefficient at $\bt^\alpha$ of $T_{\Delta+\Omega+\frac{\dd}{\dd \bt}}^{\ev}\left(\sum_{\xi}a_\xi\bt^\xi\right)$ is given by \begin{align*} &\frac{1}{\alpha!}\; \ev\left[(\delta+\partial+\frac{\dd}{\dd \bt})^\alpha(\sum_\xi a_\xi \bt^\xi)\right]=\\ &\qquad=\frac{1}{\alpha!}\ev\left[\sum_{\beta+\gamma\leq \alpha}\binom{\alpha}{\beta+\gamma}\binom{\beta+\gamma}{\beta}\delta^\beta\partial^\gamma\frac{\dd^{\alpha-\beta-\gamma}}{\dd \bt}(\sum_\xi a_\xi \bt^\xi)\right] \\ &\qquad= \frac{1}{\alpha!}\ev\left[ \sum_{\beta+\gamma\leq \alpha} \sum_\xi \binom{\alpha}{\beta+\gamma}\binom{\beta+\gamma}{\beta}\delta^\beta\partial\gamma(a_\xi)\frac{\dd^{\alpha-\beta-\gamma}}{\dd \bt}(\bt^\xi)\right] \\ &\qquad= \frac{1}{\alpha!}\sum_{\beta+\gamma\leq \alpha}\binom{\alpha}{\beta+\gamma}\binom{\beta+\gamma}{\beta}\delta^\beta\partial^\gamma(a_{\alpha-\beta-\gamma})\cdot (\alpha-\beta-\gamma)! \\ &\qquad=\sum_{\beta+\gamma\leq \alpha} \frac{1}{\beta!\cdot \gamma!} \delta^\beta\partial^\gamma(a_{\alpha-\beta-\gamma}), \end{align*} which is the term in \eqref{bothsides}. We now compute the right-hand-side of \eqref{equaltaylor}, when applied to $f$. The coefficient at $\bt^\alpha$ is \begin{align*} &\frac{1}{\alpha!}\; \ev\left[(\delta+\frac{\dd}{\dd \bt})^\alpha(T_{\Omega+\frac{\dd}{\dd \bt}}^{\ev}(\sum_{\xi}a_\xi{\bt}^\xi))\right]=\\ &\qquad= \frac{1}{\alpha!}\; \ev\left[(\delta+\frac{\dd}{\dd \bt})^\alpha(\sum_{\zeta}\frac{1}{\zeta!}\ev((\partial+\frac{\dd}{\dd \bt})^\zeta(\sum_\xi a_\xi{\bt}^\xi)){\bt}^\zeta)\right] \\ &\qquad=\frac{1}{\alpha!}\; \ev\left[(\delta+\frac{\dd}{\dd \bt})^\alpha(\sum_{\zeta}\frac{1}{\zeta!}\ev(\sum_{\gamma\leq \zeta}\binom{\zeta}{\gamma}\partial^\gamma\frac{\dd^{\zeta-\gamma}}{\dd \bt}(\sum_\xi a_\xi{\bt}^\xi)){\bt}^\zeta)\right] \\ &\qquad=\frac{1}{\alpha!}\; \ev\left[(\delta+\frac{\dd}{\dd \bt})^\alpha(\sum_{\zeta}\frac{1}{\zeta!}(\sum_{\gamma\leq \zeta}\binom{\zeta}{\gamma}\partial^\gamma(a_{\zeta-\gamma})\cdot (\zeta-\gamma)!) {\bt}^\zeta)\right] \\ &\qquad=\frac{1}{\alpha!}\; \ev\left[(\delta+\frac{\dd}{\dd \bt})^\alpha(\sum_{\zeta}\sum_{\gamma\leq \zeta}\frac{1}{\gamma!}\partial^\gamma(a_{\zeta-\gamma}) {\bt}^\zeta)\right] \\ &\qquad=\frac{1}{\alpha!}\; \ev\left[\sum_{\zeta}\sum_{\beta\leq \alpha}\sum_{\gamma\leq \zeta}\frac{1}{\gamma!}\binom{\alpha}{\beta}\delta^{\beta}\partial^\gamma(a_{\zeta-\gamma}) \frac{\dd^{\alpha-\beta}}{\dd \bt}({\bt}^\zeta)\right] \\ &\qquad=\frac{1}{\alpha!}\; \sum_{\beta\leq \alpha}\sum_{\gamma\leq \alpha-\beta}\frac{1}{\gamma!}\binom{\alpha}{\beta}\delta^{\beta}\partial^\gamma(a_{\alpha-\beta-\gamma}) \cdot(\alpha-\beta)! \\ &\qquad=\sum_{\beta\leq \alpha}\sum_{\gamma\leq \alpha-\beta}\frac{1}{\beta!\cdot\gamma!}\delta^{\beta}\partial^\gamma(a_{\alpha-\beta-\gamma}) \\ &\qquad=\sum_{\beta+\gamma \leq \alpha}\frac{1}{\beta!\cdot\gamma!}\delta^{\beta}\partial^\gamma(a_{\alpha-\beta-\gamma}), \end{align*} which is the term in \eqref{bothsides}, as required. \end{proof} \par \noindent What will be important to us is the following consequence. \begin{FACT}{Corollary}{TevDeltaPlusDDtIsIso} For any family of commuting derivations $\Delta=\{\delta_1,\dots,\delta_m\}$ on a $\Q$-algebra $A$, the Taylor morphism of the evaluation map $\ev:A[[\bt]]\lra A$ at $0$ is an isomorphism of differential rings \[ T_{\Delta+\frac{\dd}{\dd \bt}}^{\ev}:(A[[\bt]], \Delta+\frac{\dd}{\dd \bt})\to (A[[\bt]],\frac{\dd}{\dd \bt}). \] \noindent Its compositional inverse is $T_{-\Delta+\frac{\dd}{\dd \bt}}^{\ev}$, where $-\Delta$ is the family of commuting derivations $\{-\delta_1,\dots,-\delta_m\}$. Furthermore, $T^{\ev}_{\Delta+\frac{\dd}{\dd \bt}}$ is also a differential isomorphism $$(A[[\bt]], \Delta)\to (A[[\bt]],\Delta).$$ \end{FACT} \begin{proof} We recall that $T_{\Delta+\frac{\dd}{\dd \bt}}^{\ev}$ is a differential homomorphism $(A[[\bt]], \Delta+\frac{\dd}{\dd \bt})\to (A[[\bt]],\frac{\dd}{\dd \bt})$. By Theorem \ref{propertyT}, we have $$T_{\Delta+\frac{\dd}{\dd \bt}}^{\ev}\circ T_{-\Delta+\frac{\dd}{\dd \bt}}^{\ev}=T_{\frac{\dd}{\dd \bt}}^{\ev}=T_{-\Delta+\frac{\dd}{\dd \bt}}^{\ev}\circ T_{\Delta+\frac{\dd}{\dd \bt}}^{\ev}.$$ It is easy to check that $T_{\frac{\dd}{\dd \bt}}^{\ev}$ is the identity on $A[[\bt]]$. Hence, $T_{-\Delta+\frac{\dd}{\dd \bt}}^{\ev}$ is the compositional inverse of $T_{\Delta+\frac{\dd}{\dd \bt}}^{\ev}$. \par It follows that $T_{-\Delta+\frac{\dd}{\dd \bt}}^{\ev}$ is also a differential isomorphism $(A[[\bt]],\frac{\dd}{\dd \bt})\to (A[[\bt]],\Delta+\frac{\dd}{\dd \bt})$, in other words that $T_{-\Delta+\frac{\dd}{\dd \bt}}^{\ev}\circ \frac{\dd}{\dd t_i}=(\delta_i+\frac{\dd}{\dd t_i})\circ T_{-\Delta+\frac{\dd}{\dd \bt}}^{\ev}$. Now $\frac{\dd}{\dd t_i}\circ T_{-\Delta+\frac{\dd}{\dd \bt}}^{\ev}=T_{-\Delta+\frac{\dd}{\dd \bt}}^{\ev}\circ(-\delta_i+\frac{\dd}{\dd t_i})$, because $T_{-\Delta+\frac{\dd}{\dd \bt}}^{\ev}$ is a differential isomorphism $(A[[\bt]],-\Delta+\frac{\dd}{\dd \bt})\to (A[[\bt]],\frac{\dd}{\dd \bt})$. It follows that \[ T_{-\Delta+\frac{\dd}{\dd \bt}}^{\ev}\circ \frac{\dd}{\dd t_i}=\delta_i\circ T_{-\Delta+\frac{\dd}{\dd \bt}}^{\ev}+T_{-\Delta+\frac{\dd}{\dd \bt}}^{\ev}\circ(-\delta_i+\frac{\dd}{\dd t_i}), \] which implies $T_{-\Delta+\frac{\dd}{\dd \bt}}^{\ev}\circ\delta_i=\delta_i\circ T_{-\Delta+\frac{\dd}{\dd \bt}}^{\ev}$, as claimed in the ``furthermore" part. \end{proof} \par \noindent We now use \ref{TevDeltaPlusDDtIsIso} to introduce a twisting of the Taylor morphism. \begin{FACT}{+The twisted Taylor morphism.}{twistedTaylor} We assume all derivations commute. Let $A$ be a differential ring with derivations $\Delta=\{\delta_1,\dots,\delta_m\}$ and let $B$ be a $\Q$-algebra and a differential ring with derivations $\Omega=\{\partial_1,\dots,\partial_m\}$. Let $\phi:A\lra B$ be a (not necessarily differential) ring homomorphism. We write $\partial_i$ again for the extension of $\partial_i$ to $B[[\bt]]$, $\bt=(t_1,\ldots,t_m)$ obtained from differentiating coefficients as explained in \ref{SetupTaylor}. Let $\ev:B[[\bt]]\lra B$ be the evaluation map at $0$. If we equip $B[[\bt]]$ with the derivations $\Omega +\frac{\dd}{\dd \bt}$ as in \ref{SetupTaylor} and apply \ref{TevDeltaPlusDDtIsIso} for $(B,\Omega)$, we get a differential ring isomorphism\looseness=-1 \[ T^\ev_{\Omega+\frac{\dd}{\dd \bt}}:(B[[\bt]],\Omega+\frac{\dd}{\dd \bt})\lra (B[[\bt]],\frac{\dd}{\dd \bt}) \] with compositional inverse $T^\ev_{-\Omega+\frac{\dd}{\dd \bt}}$. Consequently, the map \[ T^*_\phi:=T^\ev_{-\Omega+\frac{\dd}{\dd \bt}}\circ T^\phi_\Delta:(A,\Delta)\xrightarrow{T^\phi_\Delta} (B[[\bt]],\frac{\dd}{\dd \bt})\xrightarrow{T^\ev_{-\Omega+\frac{\dd}{\dd \bt}}} (B[[\bt]],\Omega +\frac{\dd}{\dd \bt}) \] is a differential ring homomorphism $(A,\Delta)\lra (B[[\bt]],\Omega+\frac{\dd }{\dd \bt})$, called the \notion[]{twisted Taylor morphism} of $\phi$. Writing $T^*_\phi(a)=\sum_\alpha b_\alpha \bt^\alpha$, the $b_\alpha$'s are explicitly computed as \[ b_\alpha=\frac{1}{\alpha!}\sum_{\beta\leq \alpha} (-1)^{\alpha-\beta}{\alpha\choose \beta}\partial^{\alpha-\beta}\bigl(\phi(\delta^\beta(a))\bigr). \] If $a\in A$ and $\Z\{a\}$ denotes the differential subring generated by $a$ in $A$, one checks readily that \begin{enumerate}[(i)] \item $T^\phi_\Delta(a)=\phi(a)\iff \delta^\alpha (a)\in \ker(\phi)$ for all nonzero $\alpha\in \N^m$. \item $T_\phi^*(a)=T^\phi_\Delta(a)\iff \phi(\Z\{a\})$ is contained in the ring of $\Omega$-constants of $B$. \item $T_\phi^*(a)=\phi(a)\iff $the restriction of $\phi$ to $\Z\{a\}$ is a differential homomorphism. \end{enumerate} \par \noindent Hence by the implication $\La$ in (iii), if $R$ is a differential subring of $A$ such that the restriction $\phi|_R$ is a differential ring homomorphism $(R,\Delta|_R)\lra (B,\Omega)$, then $T_\phi^*$ extends $\phi$ and the part showing solid arrows in the following diagram is commutative: 
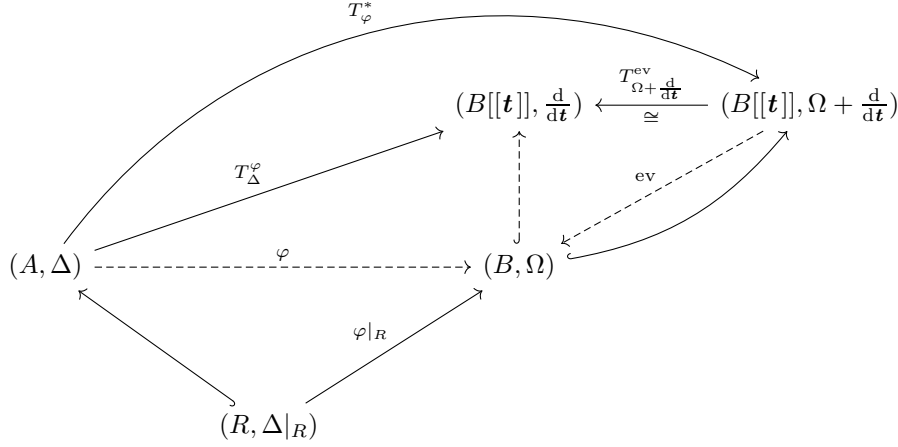
\begin{center} \begin{tikzcd}[row sep=10ex,column sep=10ex] \ & \ & (B[[\bt]],\frac{\dd}{\dd \bt}) & (B[[\bt]],\Omega+\frac{\dd}{\dd \bt})\ar[l,"T_{\Omega+\frac{\dd}{\dd \bt}}^\mathrm{ev}"', "\cong"] \ar[ld,"\mathrm{ev}"',dashed]\\ (A,\Delta)\ar[urrr, "T_\phi^*" ,bend left=40] \ar[rr,"\phi",dashed] \ar[rru,"T_\Delta^\phi"] & \ & (B,\Omega) \ar[u, hook,dashed]\ar[ur, hook, bend right=20] & \ \\ \ & (R,\Delta|_R)\ar[ul,hook]\ar[ur,"\phi|_R"] & \ & \ \end{tikzcd} \end{center} Notice that all solid arrows in this diagram are differential homomorphisms. The main case for us is when $R=K$ is a field and $B$ is a $K$-algebra such that $\phi$ is a $K$-algebra homomorphism. In this case the twisted Taylor morphism $T_\phi^*$ is in fact a differential $K$-algebra homomorphism. \par \end{FACT} \begin{FACT}{Corollary}{DiffPointFromPointInECext} Let $(K, \Delta)$ be a differential field that is large as a field and let $S$ be a differentially finitely generated $K$-algebra. If there is a $K$-algebra homomorphism $S\to L$ for some field extension $L/K$ in which $K$ is \ec\ (as a field, there are no derivations on $L$ given), then there is a differential $K$-algebra homomorphism $S\to K[[\bt]]$, where the derivations on $K[[\bt]]$ are $\Delta+\frac{\dd}{\dd \bt}$ as described above. \end{FACT} \begin{proof} Since $K$ as a field is \ec\ in $L$, there is a field extension $L'$ of $L$ which is an elementary extension of the field $K$. We replace $L$ by $L'$ if necessary and assume that $L$ is an elementary extension of the field underlying $K$. As $K$ is large, also $L$ is large. We equip $L$ with a set of commuting derivations extending those on $K$ (this is chosen arbitrarily and can always be done). \par By Theorem \ref{twistedTaylor}, there is a differential $K$-algebra homomorphism $S\to L((\bt))$. As $L$ is large and also an elementary extension of the field $K$, we know that $K$ is \ec\ as a field in $L((\bt))$. Hence by \ref{StructureThmConsequence}(ii) there is a $K$-algebra homomorphism $S\lra K$. By \ref{twistedTaylor} there is a differential $K$-algebra homomorphism $S\lra K[[\bt]]$. \end{proof} \section{Differentially large fields and algebraic characterisations} \labelx{sectionUCalgebraically} \noindent In this section we introduce the notion of differential largeness and characterize it in multiple ways, see \ref{CharDlargeI} and \ref{DLisUC}. First we recall the notion of largeness of fields. \begin{FACT}{Definition}{} A field $K$ is said to be \notion{large} (or \textit{ample} in \cite[Rem. 16.12.3]{FriJar2008}) if every irreducible affine algebraic variety $V$ over $K$ with a smooth $K$-point has a Zariski-dense set of $K$-points (equivalently, $K$ is e.c. in the function field $K(V))$. \end{FACT} \noindent Another equivalent formulation of largeness is that $K$ is \ec\ in the formal Laurent series field $K((t))$. Examples of large fields are pseudo algebraically closed fields, pseudo real closed fields and pseudo p-adically closed fields. By \cite{Pop2010} the fraction field of any Henselian local ring is large, in particular for every field $K$ and all $n\geq 1$, the power series field $K((t_1,\ldots,t_n))$ is large. \par \smallskip\noindent {\bf Convention.} Recall that for us a differential field always means a differential field in $m$ commuting derivations $\Delta=\{\delta_1,\dots,\delta_m\}$ and of characteristic 0. For a differential field $(K,\Delta)$, we equip the Laurent series field $K((\bt))$ with the natural derivations extending those on $K$; namely, $\Delta+\frac{\dd}{\dd \bt}$ as described in the previous section.\looseness=-1 \begin{FACT}{Definition}{} A differential field $K$ is said to be \notion{differentially large} if it is large as a pure field and for every differential field extension $L$ of $K$ the following implication holds: \begin{quote} If $K$ is \ec\ in $L$ as a field, then $K$ is \ec\ in $L$ as a differential field. \end{quote} \end{FACT} \medskip\noindent We now provide several algebraic characterisations of differential largeness. These characterizations resemble to some extent the characterizations of largeness of a field and serve as justification for the terminology ``differentially large". A further characterization will be given in \ref{DLisUC}. \begin{FACT}{:Characterizations of differential largeness}{CharDlargeI} Let $K=(K,\Delta)$ be a differential field. The following conditions are equivalent. \begin{enumerate}[(i),itemsep=1ex] \item\labelx{IDlarge} $K$ is differentially large. \item\labelx{IexClosedInSeries} $K$ is \ec\ in $K((\bt))$ as a differential field, where the derivations on $K((\bt))$ are the natural ones extending those on $K$. \item\labelx{IexClosedInSeriesMultiple} $K$ is \ec\ in $K((\bt_1))\ldots((\bt_k))$ as a differential field for every $k\geq 1$. \item\labelx{IPointImpliesDiffPoint} $K$ is large as a field and every differentially finitely generated $K$-algebra that has a $K$-rational point, also has a differential $K$-rational point. \item\labelx{IStructureThm} $K$ is large and every composite $K$-algebra in which $K$ is \ec\ as a field, has a differential $K$-rational point. \item\labelx{ICompositeInSeries} Every composite differential $K$-subalgebra $S$ of $K((\bt))$ has a differential $K$-rational point. \item\labelx{IComposite} $K$ is large as a field and for every composite $K$-algebra $S=A\otimes_KP$, if $A$ has a $K$-rational point, then $S$ has a differential $K$-rational point. \item\labelx{ICompositeSmoothSingle} $K$ is large as a field and for every composite $K$-algebra $S=A\otimes_KP$, if the variety defined by $A$ is smooth and if $A$ has a $K$-rational point $A\lra K$, then $S$ has a differential $K$-rational point. \item\labelx{ICompositeSmoothMany} $K$ is large as a field and for every composite $K$-algebra $S=A\otimes_KP$, if $A$ has a smooth $K$-rational point, then $S$ has a Kolchin dense set of differential $K$-rational points (cf. \ref{ECbasic}\ref{ECbasicECasKdense}(b)). \item\labelx{Ivarieties} $K$ is large as a field and for every irreducible differential variety $V$ over $K$ such that for infinitely many $r\geq 0$ the algebraic variety $\jet_r(V)$ has a smooth $K$-point, then the set of differential $K$-rational points of $V$ is Kolchin dense in $V$; in other words, for every proper closed differential subvariety $W\subseteq V$ there is a differential $K$-point in $V\setminus W$. \end{enumerate} \end{FACT} \begin{proof} \ref{IDlarge}$\Ra $\ref{IexClosedInSeriesMultiple} In the tower $K\subseteq K((\bt_1))\subseteq K((\bt_1))((\bt_2))\subseteq \ldots\subseteq K((\bt_1))\ldots((\bt_k))$ all fields are large and therefore $K$ is \ec\ in $K((\bt_1))\ldots((\bt_k))$ as a field. So by definition of differential largeness, $K$ is \ec\ in $K((\bt_1))\ldots((\bt_k))$ as a differential field. \par \par \smallskip\noindent \ref{IexClosedInSeriesMultiple}$\Ra $\ref{IexClosedInSeries} is trivial. \par \smallskip\noindent \ref{IexClosedInSeries}$\Ra $\ref{IPointImpliesDiffPoint} Since $K$ is \ec\ in $K((\bt))$ as a differential field it is also \ec\ in $K((\bt))$ as a field and so $K$ is large as a field. Let $S$ be a differentially finitely generated $K$-algebra and assume there is a point $S\lra K$. Then by \ref{twistedTaylor}, there is a differential $K$-algebra homomorphism $S\lra K[[\bt]]$. By \ref{ECbasic}\ref{ECbasicFieldField} applied to $K\subseteq K((\bt))$, assumption \ref{IexClosedInSeries} entails a differential $K$-algebra homomorphism $S\lra K$. \par \par \smallskip\noindent \ref{IPointImpliesDiffPoint}$\Ra $\ref{IStructureThm} Take $A,P$ for $S$ as in \ref{defnComposite}. Since $K$ is also \ec\ in $A$ as a field and $A$ is a finitely generated $K$-algebra, there is a $K$-algebra homomorphism $g:A\lra K$. Since $S\cong _K A\otimes _K P$ and $P$ is a polynomial $K$-algebra, $g$ can be extended to a $K$-algebra homomorphism $S\lra K$. Hence \ref{IPointImpliesDiffPoint} applies. \par \smallskip\noindent \ref{IStructureThm}$\Ra $\ref{IDlarge} Let $L$ be a differential field extension of $K$ and suppose $K$ is \ec\ in $L$ as a field. Let $S$ be a differentially finitely generated $K$-algebra, which has a differential point $f:S\lra L$. By \ref{ECbasic}\ref{ECbasicFieldField} it suffices to find a differential point $S\lra K$. By \ref{StructureThmConsequence}(i) we may replace $S$ by a composite $K$-algebra contained in $L$ and assume that $f$ is the inclusion map $S\into L$. Now \ref{IStructureThm} applies. \par \medskip\noindent Hence we know that conditions \ref{IDlarge}--\ref{IStructureThm} are equivalent. \par \par \smallskip\noindent \ref{IPointImpliesDiffPoint}$\Ra $\ref{IComposite} If $S=A\otimes_KP$ is composite and $A$ has a $K$-rational point, then as $P$ is a polynomial $K$-algebra we may extend this point to a point $S\lra K$. By \ref{IPointImpliesDiffPoint}, $S$ has a differential $K$-rational point. \par \smallskip\noindent \ref{IComposite}$\Ra $\ref{ICompositeInSeries} If $S=A\otimes_KP$ is a composite $K$-subalgebra of $K((\bt))$, then as $K$ is a large field, $K$ is \ec\ in $A$ as a field and thus $A$ has a $K$-rational point. Now \ref{IComposite} applies. \par \smallskip\noindent \ref{ICompositeInSeries}$\Ra $\ref{IexClosedInSeries} follows from \ref{StructureThmConsequence}(i) using the characterization \ref{ECbasic}\ref{ECbasicFieldField} of \ec \par \medskip\noindent Hence we know that conditions \ref{IDlarge}--\ref{IComposite} are equivalent. \par \par \smallskip\noindent \ref{IDlarge}$\Ra $\ref{ICompositeSmoothMany} If $S=A\otimes_KP$ is composite and $A$ has a smooth $K$-rational point, then as a large field, $K$ is \ec\ in $A$ as a field. Since $P$ is a polynomial $K$-algebra we know that $S$ is a polynomial $A$-algebra and so $A$ is \ec\ in $S$ as a ring. It follows that $K$ is \ec\ in $S$ as a field and by \ref{IDlarge} (invoke \ref{ECbasic}\ref{ECbasicInDomainToField}) it is then also \ec\ in $S$ as a differential field. By \ref{ECbasic}\ref{ECbasicECasKdense} we see that $S$ has a Kolchin dense set of differential $K$-rational points. \par \smallskip\noindent \ref{ICompositeSmoothMany}$\Ra $\ref{ICompositeSmoothSingle} is trivial. \par \smallskip\noindent \ref{ICompositeSmoothSingle}$\Ra $\ref{IStructureThm} Let $S=A\otimes_KP$ be a composite $K$-algebra in which $K$ is \ec\ as a field. Then $K$ is also \ec\ in $A$ as a field and therefore it possesses a smooth $K$-rational point $f:A\lra K$ (cf. \ref{ECbasic}\ref{ECbasicECsmooth}). Pick $h\in A$ with $f(h)\neq 0$ such that the variety defined by the localization $A_h$ is smooth. We may now apply \ref{ICompositeSmoothSingle} to the composite algebra $S_h=A_h\otimes_KP$. \par \medskip\noindent Hence we know that \ref{IDlarge}--\ref{ICompositeSmoothMany} are equivalent. Property \ref{Ivarieties} is just a reformulation of the definition of differential largeness in geometric form as follows. Let $S=K\{x\}/\Dp$, $x=(x_1,\ldots,x_n)$ be a differentially finitely generated $K$-algebra and a domain with residue map $\pi:K\{x\}\lra S$. Let $V$ be the differential variety defined by $S$, hence $V=\{a\in M^n\st \Dp(a)=0\}$, where $M$ is the differential closure of $K$. Then $V$ is a $K$-irreducible differential variety defined over $K$. Now for $r\in \N$, the variety $\jet_r(V)$ has coordinate ring $A_r:=\pi(K\{x\}_{\leq r})$\footnote{Here $K\{x\}_{\leq r}$ denotes the subring of $K\{x\}$ of all polynomials in $\theta x_i$, where $\ord(\theta)\leq r$.} and $S$ is the union of the chain $(A_r)_r$ of $K$-subalgebras of $S$. Clearly $K$ is \ec\ in $S$ as a field if and only if $K$ is \ec\ in $A_r$ as a field for all (or infinitely many) $r$. Since $K$ is large, this is equivalent to saying that $\jet_r(V)$ has a smooth $K$-point for all (or infinitely many) $r$. Hence the assumption about $V$ in \ref{Ivarieties} precisely says that $K$ is \ec\ in $S$ as a field. \par On the other hand, the conclusion about $V$ in \ref{Ivarieties} precisely says that $K$ is \ec\ in $S$ as a differential field (use \ref{ECbasic}\ref{ECbasicECasKdense}). \par This shows that differential largeness is equivalent to \ref{Ivarieties} formulated for affine differential varieties. But obviously the affine case implies \ref{Ivarieties} in full. \end{proof} \begin{FACT}{Corollary}{DropDerivations} If $K=(K,\delta_1,\ldots,\delta_m,\partial_1,\ldots,\partial_k)$ is a differentially large field, $m,k\geq 0$, then also $K=(K,\delta_1,\ldots,\delta_m)$ is differentially large. \end{FACT} \begin{proof} This is immediate from the power series characterization \ref{CharDlargeI}\ref{IexClosedInSeries}. \end{proof} \begin{FACT}{Corollary}{ConstantsLargeANDec} Let $K=(K,\delta_1,\ldots,\delta_m,\partial_1,\ldots,\partial_k)$ be differentially large, $m\geq 0,k\geq 1$ and let $C=\{a\in K\st \partial_1(a)=\ldots=\partial_k(a)=0\}$ be the constant field of $(\partial_1,\ldots,\partial_k)$. \begin{enumerate}[(i)] \item $C$ is closed under the derivations $\delta_1,\ldots,\delta_m$ and $(C,\delta_1,\ldots,\delta_m)$ is \ec\ in $(K,\delta_1,\ldots,\delta_m)$. \item $(C,\delta_1,\ldots,\delta_m)$ is differentially large; when $m=0$, this just says that $C$ is a large field. \end{enumerate} \end{FACT} \begin{proof} We write $\delta=(\delta_1,\ldots,\delta_m)$ and by a trivial induction we may assume that $k=1$. Set $\partial=\partial_1$. \par \smallskip\noindent (i) Since all derivations commute, $C$ is closed under all derivations. Let $(S,\hat \delta)$ be a $(C,\delta)$-algebra that is finitely generated as such. Suppose we are given a differential $K$-rational point $\lambda:(S,\hat\delta)\lra (K,\delta)$ (in fact we will only need that $S$ has a $K$-rational point). It suffices to find a differential $C$-algebra homomorphism $(S,\hat\delta)\lra (C,\delta)$. We expand $(S,\hat\delta)$ by the trivial derivation and obtain a differentially finitely generated $(C,\delta,\partial )$-algebra $(S,\hat\delta,0)$ (note that $\partial $ is trivial on $C$). \par A straightforward calculation shows that $(S,\hat\delta,0)\otimes _C(K,\delta,\partial)$ is a differential $(K,\delta,\partial)$-algebra (the derivations are given by $\hat\delta_i\otimes \delta_i$ and $0\otimes \partial$) that is finitely generated as such, and $\lambda\otimes \id :S\otimes _CK\lra K$ is a (not necessarily differential) $K$-algebra homomorphism; also see \cite[\S3.1]{LeSTre2020} for generalities on derivations and tensor products. \par Since $(K,\delta,\partial)$ is differentially large, \ref{CharDlargeI}\ref{IPointImpliesDiffPoint} gives a differential point $\mu:(S,\hat\delta,0)\otimes _C(K,\delta,\partial)\lra (K,\delta,\partial)$ and we get a $C$-algebra homomorphism $\mu_0:S\lra S\otimes _CK\xrightarrow{\ \mu\ }K$. Since the natural map $S\lra S\otimes _CK$ is differential for $\delta$, also $\mu_0$ is a differential homomorphism $(S,\hat\delta)\lra (K,\delta)$. But $\mu_0$ has values in $C$ because for $s\in S$ we have $\partial(\mu_0(s))=\partial(\mu(s\otimes 1))=\mu((0\otimes \partial) (s\otimes 1))=\mu(0)=0$. Hence indeed $\mu_0(s)\in C$ as required. \par \smallskip\noindent (ii) Since $(K,\delta,\partial)$ is differentially large, it is \ec\ in $K((t_1,\ldots,t_{m+1}))$ when the latter is furnished with the natural derivations, cf. \ref{CharDlargeI}\ref{IexClosedInSeries}. By (i), $(C,\delta)$ is \ec\ in $(K,\delta)$. If $m=0$ it follows that $C$ is \ec\ as a field in $K((t_1))$, hence $C$ is a large field. If $m\geq 1$, we see that $(C,\delta)$ is \ec\ in $C((t_1,\ldots,t_m))$, which shows that it is differentially large by \ref{CharDlargeI}\ref{IexClosedInSeries}. \end{proof} \par \newpage \noindent At the end of this section we show that differentially large fields are first-order axiomatizable; in other words, the class of differentially large fields is an elementary class in the language of differential rings. We show this implicitly in \ref{DLisUC}, by proving that differentially large fields are precisely those large and differential fields satisfying the axiom scheme UC in \cite[4.5]{Tressl2005}; thus we refer to this paper for explicit axioms. The proof of \ref{DLisUC} only uses properties of models of UC and results from this paper. \begin{FACT}{Remark}{elementaryclass} It is worth mentioning (for the non-logician) the benefits of knowing that a class of structures is elementary (i.e., first-order axiomatisable). In our context this means that two properties hold: (1) ultraproducts of differentially large fields are again differentially large, and (2) differential fields that are existentially closed in some differentially large field are themselves differentially large. Property (2) is obvious from the characterization \ref{CharDlargeI}\ref{IexClosedInSeries}. So it is only property (1) that needs to be established. Being an elementary class, opens up the model theoretic toolbox to the analysis of differentially large fields, and it implies for example the following transfer principle (phrased in technical terms in \ref{UCrecall} below). \begin{quote} If $K$ is a differentially large field and $K$ as a pure field has ``good" elimination theory, then the differential field $K$ also has good elimination theory. \end{quote} To illustrate what ``good" elimination theory means, we look at classical examples of ``good" elementary classes of fields. Algebraically closed fields have good elimination theory, this is due to Chevalley's theorem which says that the projection of a variety is constructible. If $K$ is a real closed field or a p-adically closed field, then projections of $K$-varieties (by which we mean here Zariski closed subset of some $K^n$) are generally not constructible; however, the following weaker statement holds: the complement of a projection of a $K$-variety is again the projection of a $K$-variety (this property of a field is called ``model-completeness", cf. \cite[sectoin 8.3]{Hodges1993}). So then the transfer principle above says that for a differentially large field $K$ the following holds: if $K$ is algebraically closed as a field, then the projection of a differential variety is differentially constructible (i.e., a finite Boolean combination of Kolchin closed sets); if $K$ is real closed or p-adically closed, then the complement of a projection of a differential variety is again the projection of a differential variety. \end{FACT} \begin{FACT}{Proposition}{DLisUC} Let $K$ be a differential field that is large as a field. Then $K$ is differentially large if and only if it satisfies the axiom scheme $\UC$ from \cite[4.5]{Tressl2005}. \end{FACT} \begin{proof} First assume that $K$ is differentially large. By \cite[Theorem 6.2(II)]{Tressl2005}, there is a differential field extension $L$ of $K$ such that $L\models \UC$ and such that $K$ is elementary in $L$ as a field. In particular $K$ is \ec\ in $L$ as a field. Since $K$ is differentially large, $K$ is \ec\ in $L$ as a differential field. By \cite[Proposition 6.3]{Tressl2005}, $\UC$ has an inductive axiom system in the language of differential rings. But then $K$ also satisfies these axioms. Hence $K\models \UC$. \par \smallskip For the converse assume that $K$ is a model of $\UC$. We verify the definition of differentially large. Let $L$ be a differential field extension of $K$ such that $K$ is \ec\ in $L$ as a field. Then there is a field $M$ extending $L$ such that $K$ is elementary in $M$ as a field. In particular $M$ is a large field. We may now extend the derivations of $L$ arbitrarily to commuting derivations of $M$. Hence we may replace $L$ by $M$ furnished with these derivations and assume that $L$ is large as a field. By \cite[Theorem 6.2(II)]{Tressl2005} again, there is a differential field extension $F$ of $L$ such that $F\models \UC$ and such that $L$ is elementary in $F$ as a field. Then $K$ is \ec\ in $F$ as a field and $K,F\models \UC$. By \cite[Theorem 6.2(I)]{Tressl2005}, this shows that $K$ is \ec\ in $F$ (as a differential field), showing the assertion. \end{proof} \par \noindent By \ref{DLisUC} we may now record important properties of differentially large fields (that follow from being models of $\UC$, see \cite{Tressl2005}). \begin{FACT}{Corollary}{UCrecall} \begin{enumerate}[(i)] \item If $L$ and $M$ are differentially large fields and $K$ is a common differential subfield of $K$ then $L$ and $M$ have the same existential theory over $K$ (meaning they solve the same systems of differential equations with coefficients in $K$) if and only if they have the same existential theory over $K$ as fields. \item If $K$ is a differential field that is large as a field, then there is a differential field extension $L$ of $K$ such that $L$ is differentially large and an elementary extension of $K$ as a field. \item Let $K$ be a differentially large field and let $A\subseteq K$. Suppose $K$ is model complete as a field in the language $\SL_\mathrm{ri}(A)$ of rings extended by constant symbols naming the elements of $A$. \par Then also $K$ is model complete in the language $\SL_\mathrm{diff}(A)$ of differential rings extended by all constant symbols naming the elements of $A$. Furthermore, if $\hat \SL$ is a language extending $\SL_\mathrm{ri}$ and $\hat K$ is an expansion of $K$ to $\hat \SL_\mathrm{diff}$ such that the new symbols are $A$-definable in the field $K$ and such that the restriction of $\hat K$ to $\SL^*(A)$ has quantifier elimination \footnote{An example of $\hat \SL$ is the language $\SL_\mathrm{ri}(\leq)$ of ordered rings, $K=(K,\delta)$ is a real closed field furnished with commuting derivations, $A=\0$ and $\hat K=(K,\leq,\delta)$. The restriction of $\hat K$ to $\SL_\mathrm{ri}(A)$ then is the ordered field $(K,\leq)$.}, then $\hat K$ has quantifier elimination in the language $\hat \SL_\mathrm{diff}(A)$. \end{enumerate} \end{FACT} \section{Fundamental properties and constructions}\label{further} \noindent {\bf Summary} We show that algebraic extensions of differentially large fields are again differentially large by invoking differential Weil descent in \ref{algextdifflarge}. Specifically differentially closed fields are identified as precisely the algebraic closures of differentially large fields; in a similar way, M. Singer's closed ordered differential fields are characterized, see \ref{DLisACFRCF}. We show that a differentially large field is pseudo algebraically closed just if it is pseudo differentially closed, see \ref{psfDCF}. We characterize the existential theory of differentially large fields in \ref{existentialresult}. We show that differentially large fields are Picard-Vessiot closed in \ref{PVclosed}. In \ref{KolDenseInDiffGroups} we establish Kolchin-denseness of rational points in differential algebraic groups. \par We start with a concrete method to construct differentially large fields. This is deployed in \ref{ConstructionPowerSeriesCont} to obtain concrete constructions of differentially closed fields. \begin{FACT}{Proposition}{ConstructionDirectLimit} Let $(K_i,f_{ij})_{i,j\in I}$ be a directed system of differential fields and differential embeddings with the following properties. \begin{enumerate}[(a)] \item\labelx{ConPropLarge} All $K_i$ are large as fields. \item\labelx{ConPropEC} All embeddings $f_{ij}:K_i\lra K_j$ are isomorphisms onto a subfield of $K_j$ that is \ec\ in $K_j$ as a field. \item\labelx{ConPropSeries} If $i\in I$, then there is some $j\geq i$ such that a differential homomorphism $K_i[[\bt]]\lra K_j$ extending $f_{ij}$ exists. \end{enumerate} Then the direct limit $L$ of the directed system is a differentially large field. \end{FACT} \begin{proof} We write $f_i:K_i\lra L$ for the natural map into the limit, which obviously is a differential homomorphism between differential fields. We use the characterization \ref{CharDlargeI}\ref{IComposite} to show that $L$ is differentially large. Firstly, $L$ is large as a field, because if $C$ is a curve defined over $L$ that has a smooth $L$-rational point then take $i\in I$ such that $C$ is defined over $K_i$ (via $f_i$) and such that $C$ has a smooth $K_i$-rational point. By \ref{ConPropLarge}, the curve $C$ has infinitely many $K_i$-rational points and so it also has infinitely many $L$-rational points. \par Now let $S$ be a differentially finitely generated $L$-algebra and a domain that has a point $S\lra L$. Pick $r\in \N$ and a differential prime ideal $\Dp$ of $L\{x\}$, $x=(x_1,\ldots,x_r)$ such that $S=L\{x\}/\Dp$. By the Ritt-Raudenbusch basis theorem there is a finite $\Sigma\subseteq \Dp$ whose differential radical is $\Dp$. By \ref{CharDlargeI}\ref{IComposite} it suffices to find a differential zero of $\Sigma$ in $L$. Take $i\in I$ with $\Sigma \subseteq f_i(K_i)\{x\}$ and let $S_0:=K_i\{x\}/f_i^{-1}(\Dp)$. Then $S_0$ is a differentially finitely generated $K_i$-algebra and the composition of the natural embedding $S_0\lra S$ with a point $S\lra L$ is a homomorphism $S_0\lra L$ extending $f_i$. We now want to invoke \ref{DiffPointFromPointInECext} and here we need \ref{ConPropEC}. Namely, with this condition one readily verifies definition \ref{ECbasic} and checks that $K_i$ is existentially closed in $L$ as a field (via $f_i$). \par Hence we may apply \ref{DiffPointFromPointInECext} to obtain a differential $K$-algebra homomorphism $S_0\lra K_i[[\bt]]$. Finally assumption \ref{ConPropSeries} gives us a differential $K_i$-algebra homomorphism $S_0\lra K_j$ for some $j\geq i$. This yields a differential solution of $\Sigma $ in $L$. \end{proof} \par \noindent Concretely, \ref{ConstructionDirectLimit} may be used to produce differentially large fields via iterated power series constructions using standard power series, Puiseux series or generalised power series. Here are two instances. \begin{FACT}{Examples}{ConstructionPowerSeries} Let $K$ be a differential field. We write $K_0=K$. \begin{enumerate}[(i),itemsep=1ex] \item\labelx{ConstructionPowerSeriesA} We define by induction on $n\geq 0$, the differential field extension $K_{n+1}$ of $K_n$ as $K_{n+1}=K_n((\bt_n))$, where $\bt_n=(t_{n1},\ldots,t_{nm})$; the derivations on $K_{n+1}$ are the natural ones, extending those on $K_n$ and satisfying $\delta_j(t_{nk})=\frac{\dd}{\dd t_{nj}}(t_{nk})$. Then $L=\bigcup _{n\in\N}K_n$ is differentially large. If $K$ is large as a field, then $K$ is \ec\ in $L$ as a field. \par To see this we apply \ref{ConstructionDirectLimit} to the family of all $K_n$, $n>0$ together with the inclusion maps $K_i\into K_j$ for $i\leq j$. Hence $L$ is differentially large. Since all $K_n$ are large fields we know that they are \ec\ in $L$ as a field. Hence if $K$ happens to be large as a field, then $K$ is also \ec\ in $L$ as a field. \par \smallskip\noindent This example is continued in \ref{ConstructionPowerSeriesCont}. \item\labelx{ConstructionPowerSeriesB} Assume here that the number $m$ of derivations is 1. Then the generalised power series field $K((t^\Q))$ carries a derivation defined by $\frac{\dd}{\dd t}(\sum a_\gamma t^\gamma)=\sum{a_\gamma\mal \gamma}\mal t^{\gamma-1}$ and the given derivation $\delta$ on $K$ can be extended to a derivation $\partial$ by $\partial(\sum a_\gamma t^\gamma)=\sum \delta(a_\gamma) t^\gamma$. We consider $K((t^\Q))$ as a differential field extension of $K((t))$, equipped with the derivation $\frac{\dd}{\dd t}+\partial$. \par Now define $K_{n+1}=K_n((t_n^\Q))$. Since $K_n$ carries a henselian valuation for $n>0$ we know that $K_n$ is a large field. Hence conditions \ref{ConstructionDirectLimit}\ref{ConPropLarge},\ref{ConPropSeries} hold for the family of all $K_n$, $n>0$ and the inclusion maps $K_i\into K_j$ when $i\leq j$. \par If $K$ is algebraically closed, real closed or p-adically closed, then so are all $K_n$ and by standard theorems from model theory condition \ref{ConstructionDirectLimit}\ref{ConPropEC} holds in each case. Thus $K_\infty=\bigcup_nK_n$ is a differentially large field. Furthermore $K_\infty=\bigcup_nK_n$ is again algebraically closed, real closed or p-adically closed, respectively. To be precise: If $K$ is algebraically closed, then $K_\infty$ is a differentially closed field; if $K$ is real closed, then $K_\infty$ is a closed ordered differential field in the sense of \cite{Singer1978a}; and if $K$ is p-adically closed, then $K_\infty$ is an existentially closed differential field in the class of p-adically valued and differential fields as considered in \cite{GuzPoi2010}. \end{enumerate} \end{FACT} \begin{FACT}{Remark}{} In Example~\ref{ConstructionPowerSeries}(ii) one can replace generalised power series with Puiseux series. More precisely, if $K_{n+1}$ is now defined to be the Puiseux series field over $K_n$, namely $$K_{n+1}=K_n((t^{\frac{1}{\infty}}))=\bigcup_{k\in\N}K_n((t^{\frac{1}{k}})),$$ then exactly the same conclusions can be derived on $K_\infty$. \end{FACT} \begin{FACT}{Counterexample}{IterateAlgPowerSeries} One cannot replace iterated power series by iterated algebraic power series in \ref{ConstructionPowerSeries}\ref{ConstructionPowerSeriesA}: Let $K$ be an ordinary differential field and let $L=\bigcup _nK((t_1))_\alg\ldots ((t_n))_\alg$, where the derivation is chosen as in \ref{ConstructionPowerSeries}\ref{ConstructionPowerSeriesA}. Thus $L$ is a differential subfield of $K(t_1,t_2,\ldots)^\alg$, where $\delta(t_i)=1$ for all $i$. Notice that $L$ is a large field, because the natural valuation of algebraic power series is henselian. If $L$ were differentially large, then $L$ has non-trivial solutions of the differential equation $\delta x=x$. However, we show that this is in general not the case even for $M=K(t_1,t_2,\ldots)^\alg$. \par To see this, consider the following property of a differential large field $F$. \[ \forall x\in F, n\in \N_{>0}: \delta(x)=n\mal x\Ra x=0.\leqno(\dagger) \] Then, if $F$ has property $(\dagger)$ so does its algebraic closure $F^\alg$ and its function field $F(t)$, where $\delta(t)=1$. Hence if we start with $K$ being a differential field with trivial derivation, then by induction, property $(\dagger)$ passes to $K(t_1,\ldots,t_n)^\alg$ and so also passes to $M$. \par For the proof that $(\dagger)$ passes to $F(t)$, assume that $\delta(f/g)=n\mal f/g$ with $g$ monic and $f$ with leading coefficient $a$. Then $nfg=\delta(f)g-f\delta(g)=(f^\delta+f')g-f(g^\delta+g')$ and comparing leading coefficients shows that $n\mal a=\delta(a)$. Hence by $(\dagger)$ for $F$ we get $a=0$ as required. \par For the proof that $(\dagger)$ passes to $F^\alg$, one first checks that it passes to $F(C)$, where $C$ is the constant field of $F^\alg$. Hence we may replace $F$ by $F(C)$ and assume that $F$ and $F^\alg$ have the same constant field. Let $\alpha$ be algebraic over $F$ with minimal polynomial $f$ and assume $\delta(\alpha)=n\mal \alpha$. then any other root $\beta$ of $f$ also satisfy this equation, which implies that $\delta(\frac{\alpha}{\beta})$ is a constant, thus it is in $F$. Hence $F(\alpha)$ is the splitting field of $f$ and so $F(\alpha)/F$ is Galois. Let $d$ be the order of the Galois group and let $\sigma\in\Gal(F(\alpha)/F)$. As we have seen, $\sigma(\alpha)=c\mal \alpha$ for some constant $c$. Hence $\alpha=\sigma^d(\alpha)=c^d\mal \alpha$ and so $c^d=1$. But then $\sigma(\alpha^d)=(c\alpha)^d=\alpha^d$, which shows that $\alpha^d$ is in the fixed field $F$. Since $\delta(\alpha^d)=d\mal n\mal \alpha^d$, we get $\alpha=0$ from $(\dagger)$ for $F$. \end{FACT} \begin{FACT}{:The existential theory of differentially large fields.}{DLexTheory} The existential theory of the class of all large fields of characteristic 0 is the existential theory of the field $\Q((t))$ (cf. \cite[Prop. 2.25]{Sander1996a}). This follows essentially from the fact that $\Q((t))$ is itself a large field. Since the existential theory of a differentially large field is uniquely determined by its existential theory of its field structure -- in the sense of \ref{UCrecall}(i) -- one is led to the question on whether the existential theory of the class of differentially large fields is the existential theory of $\Q((\bt ))$, equipped with its natural derivations. \par However $\Q((\bt ))$ does not satisfy the existential theory of the class of differentially large fields (and so it is not differentially large either). To see an example, let $C$ be the curve defined by $x^3+y^3=1$. Then $(1,0)$ and $(0,1)$ are the only rational points on $C$ (and they are regular points). Hence the sentence $\phi$ saying that there is a point $(x,y)$ on $C$ with $x\neq 0,y\neq 0$ and $x'=y'=0$ fails in the differential field $\Q((t))$ (we work with $m=1$ here). On the other hand $\phi$ is true in every differentially large field $K$, because the constants of $K$ are large as a field by \ref{ConstantsLargeANDec}. \end{FACT} \medskip\noindent On the positive side we now show that \begin{FACT}{Theorem}{existentialresult} The existential theory of the class of differentially large fields is the existential theory of $\Q((\bt_1))((\bt_2))$. \end{FACT} \begin{proof} Let $\Sigma\subseteq \Z\{x_1,\ldots,x_n\}$ be a system of differential polynomials in $n$ variables and $m$ commuting derivations. If $\Sigma$ has a solutions in $\Q((\bt_1))((\bt_2))$ and $K$ is a differentially field, then $\Sigma$ also has a solution in $K((\bt_1))((\bt_2))$. Hence if $K$ is differentially large, then by \ref{CharDlargeI}\ref{IexClosedInSeriesMultiple}, $\Sigma$ also has a solution in $K$. \par Conversely, suppose $\Sigma$ has a solution in every differentially large field. By \ref{UCrecall}(ii) there is a differentially large field $K$ containing $\Q((\bt_1))$ as a differential subfield such that the extension $K/\Q((\bt_1))$ of fields is elementary. Let $S_0=\Q\{x_1,\ldots,x_n\}/\sqrt[d]\Sigma$ and let $f:S_0\lra K$ be a differential point of $S_0$. Let $\Dp=\Ker(f)$ and let $S=S_0/\Dp$. It suffices to find a differential point $S\lra \Q((\bt_1))((\bt_2))$. Write $S=A_h\otimes_\Q P$ as in \ref{structuretheorem}. The the restriction $f|_{A_h}$ is a $K$-rational point of $A_h$. Since $A_h$ is a finitely generated $\Q$-algebra and $\Q((\bt_1))$ is \ec\ in $K$ as a field, there is also a point $g_0:A_h\lra \Q((\bt_1))$. Since $P$ is a polynomial $\Q$-algebra, $g_0$ can be extended to a point $g:S\lra \Q((\bt_1))$. By \ref{twistedTaylor}, there is a differential point $S\lra \Q((\bt_1))[[\bt_2]]$. \end{proof} \begin{FACT}{:Differentially large fields are PV-closed.}{PVclosed} We prove that differentially large fields solve plenty of algebraic differential equations. Namely, we prove that they solve all consistent systems of linear differential equations. We first show that they are Picard-Vessiot closed (or PV-closed). \par Let $(K,\delta_1,\dots,\delta_m)$ be a differential field, and let $A_i\in Mat_n(K)$, for $i=1,\dots,m$ satisfying what is called the \textit{integrability condition}; namely \[ \delta_iA_j-\delta_j A_i=[A_i,A_j], \] where $\delta_j A_j$ denotes the $n\times n$ matrix obtained by applying $\delta_i$ to $A_j$ entry-wise. The differential field $K$ is said to be PV-closed if for each such tuple $(A_1,\dots,A_m)$ of matrices there is a $Z\in GL_n(K)$ such that $$\delta_i Z=A_i Z_i \quad \text{ for } i=1,\dots,m.$$ \begin{FACT}{Lemma}{DiffLargeImpliesPVclosed} Every differentially large field is PV-closed. \end{FACT} \begin{proof} Suppose $K$ is a differentially large field. Suppose $A_1,\dots,A_m$ are elements in $Mat_n(K)$ satisfying the integrability condition. Let $X$ be an $n\times n$ matrix of variables and define derivations on $K(X)$ that extend the ones in $K$ and satisfy $$\delta_i X=A_i X$$ Then, by the integrabilty condition, these derivations commute in all of $K(X)$. Since $K$ is e.c. in $K(X)$ as fields, by differential largeness, it is also e.c. as differential fields. This yields the desired (fundamental) solution in $K$. \end{proof} \par \noindent In differentially large fields, \ref{DiffLargeImpliesPVclosed} is a special case of a stronger property: \begin{FACT}{Proposition}{solvelinear} Let $\Sigma$ and $\Gamma$ be finite collections of differential polynomials in $K\{x_1,\dots,x_n\}$. Assume that the system $$P=0 \; \&\; Q\neq 0, \quad \text{ for } P\in \Sigma \text{ and } Q\in \Gamma$$ is consistent (i.e., it has a solution in some differential field extension of $K$). If $\Sigma$ consists of linear differential polynomials and $K$ is differentially large, then the system has a solution in $K$. \end{FACT} \begin{proof} Since the system is assumed to be consistent, the differential ideal generated by $\Sigma$ in $K\{x_1,\dots,x_n\}$, denoted $[\Sigma]$, is prime. Thus, the differential field extension $L=$qf$(K\{x_1,\dots,x_n\}/[\Sigma])$ has a solution to the system. Since $[\Sigma]$ is generated, as an ideal of $K\{x_1,\dots,x_n\}$, by linear terms, we get that $K$ is e.c. in $L$ as fields, and, by differential largeness, also as differential fields. The result follows. \end{proof} \par \end{FACT} \begin{FACT}{:A glimpse on the Differential Weil Descent.}{WeilGlimpse} If $K$ is a large field, then every algebraic field extension of $K$ is again large. This follows from an argument involving Weil descent in the case when $L/K$ is finite, see \cite[Theorem 2.14]{Sander1996a} and \cite[Prop. 1.2]{Pop1996}. For differentially large fields, this can also be carried out. We will explain a special case of the differential Weil descent suitable for our purpose and refer to \cite[Theorem 3.4]{LeSTre2020} for the general assertion and for proofs. \par We will be working with a finite extension $L/K$ of differential fields and a differential $L$-algebra $S$. Then the classical Weil descent $W(S)$ of the underlying $L$-algebra of $S$ is a $K$-algebra and there is a ``natural" bijection \[ \Hom_{K\text{-}\kat{Alg}}(W(S),K)\lra \Hom_{L\text{-}\kat{Alg}}(S,L). \] Here homomorphisms are algebra homomorphisms over $K$ and $L$, respectively. Now in \cite[Theorem 3.4]{LeSTre2020} it is shown that the ring $W(S)$ can be naturally expanded to a differential $K$-algebra $W^{diff}(S)$ such that the bijection above restricts to a bijection \[ \Hom_{\text{diff. }K\text{-}\kat{Alg}}(W^\mathrm{diff}(S),K)\lra \Hom_{\text{diff. }L\text{-}\kat{Alg}}(S,L). \] This time, homomorphisms are differential algebra homomorphisms over $K$ and $L$, respectively. The terminology ``natural" in both bijections refers to the fact that in fact $W$ and $W^{diff}$ are functors defined on categories of algebras and differential algebras, respectively. However for our application below only the existence of the bijections above are needed. We refer to \cite[Section 3]{LeSTre2020} for a self contained exposition of the matter, where all data are constructed explicitly. In particular the construction there shows that $W^\mathrm{diff}(S)$ is a differentially finitely generated $K$-algebra if $S$ is a differentially finitely generated $L$-algebra. \end{FACT} \begin{FACT}{Theorem}{algextdifflarge} If $K$ is differentially large, then so is every algebraic extension (equipped with the induced derivations). \end{FACT} \begin{proof} Let $L/K$ be an algebraic extension. We first deal with the case when $L/K$ is finite. We verify condition \ref{CharDlargeI}\ref{IPointImpliesDiffPoint} for $L$. So let $S$ be a differentially finitely generated $L$-algebra that has an $L$-rational point. Let $W^{diff}(S)$ be the differential Weil descent as explained in \ref{WeilGlimpse}. Thus, $W^{diff}(S)$ is a differentially finitely generated $K$-algebra and we have a bijection \[ \Hom_{K\text{-}\kat{Alg}}(W(S),K)\lra \Hom_{L\text{-}\kat{Alg}}(S,L), \] which restricts to a bijection \[ \Hom_{\text{diff. }K\text{-}\kat{Alg}}(W^\mathrm{diff}(S),K)\lra \Hom_{\text{diff. }L\text{-}\kat{Alg}}(S,L). \] Since $S$ has an $L$-rational point we may use the first bijection and see that $W(S)$ has a $K$-rational point. Since $K$ is differentially large there is a differential $K$-rational point $W^\mathrm{diff}(S)\lra K$. Using the second bijection we see that $S$ has a differential $L$-rational point. \par \medskip Hence we know the assertion when $L/K$ is finite. In general, let $S=A\otimes P$ be a composite $L$-algebra such that the affine variety defined by $A$ is smooth. Suppose there is an $L$-rational point $A\lra L$. By \ref{CharDlargeI}\ref{ICompositeSmoothSingle} it suffices to show that there is a differential point $S\lra L$. Write $S=L\{x\}/\Dp$, $x=(x_1,\ldots,x_r)$, for a prime differential ideal $\Dp$ of $L\{x\}$ and let $\Sigma\subseteq \Dp$ be finite with $\Dp=\sqrt[d]\Sigma$. It suffices to find a differential solution of $\Sigma=0$ in $L$. Choose a finite extension $K_0/K$ in $L$ with $\Sigma\subseteq K_0\{x\}$. Let $S_0=K_0\{x\}/\Dp\cap K_0\{x\}$, which we consider as a subring of $S$. By \ref{structuretheorem} there are a finitely generated $K_0$-subalgebra $A_0$ of $S_0$, a polynomial $K_0$-subalgebra $P_0$ of $S_0$ and an element $h\in A_0$ such that $(S_0)_h\cong (A_0)_h\otimes_{K_0}P_0$. \par Since $A_0\subseteq S$ is finitely generated we may write $P=P_1\otimes_L P_2$ for some polynomial $L$-algebras $P_i$, $P_1$ finitely generated such that $A_0\subseteq A\otimes_L P_1$. Then $A\otimes_L P_1$ is again finitely generated, the affine variety defined by $A\otimes_L P_1$ is again smooth and still has an $L$-rational point. Since $L$ is large, there is also an $L$-rational point $(A\otimes_L P_1)_h\lra L$. Via restriction we get an $L$-rational point $f:(A_0)_h\lra L$. Since $(A_0)_h$ is finitely generated as a $K_0$-algebra, there is a finite extension $K_1/K_0$ contained in $L$ such that $f$ has values in $K_1$. Since $P_0$ is a polynomial $K_0$-algebra, $f$ can be extended to a $K_1$-rational point $(S_0)_h\lra K_1$. Tensoring with $K_1$ gives a $K_1$-rational point of $(S_0)_h\otimes_{K_0}K_1$. The latter is a differentially finitely generated $K_1$-algebra. By what we have shown, $K_1$ is differentially large. By \ref{CharDlargeI}\ref{IPointImpliesDiffPoint} there is a differential point $(S_0)_h\otimes_{K_0}K_1\lra K_1$. Since $\Sigma\subseteq K_1\{x\}$ this gives rise to a differential solution of $\Sigma=0$ in $K_1\subseteq L$. \end{proof} \begin{FACT}{Corollary}{DLisACFRCF} The algebraic closure of a differentially large field is differentially closed. In particular, if $K\models \operatorname{CODF}_m$, the theory of closed ordered differential fields in $m$ commuting derivations, then $K(i)\models \DCF_{0,m}$. \end{FACT} \noindent Previously known examples of differential fields with minimal differential closures are models of \CODF\ (which we denote as $\CODF_1$), see \cite{Singer1978b}, and fixed fields of models of $\DCF_{0,m}\operatorname{A}$, the theory differentially closed fields with a generic differential automorphism, see \cite{LeoSan2016}. The corollary delivers a vast variety of new differential fields with this property, namely all differentially large fields, see also \ref{UCrecall}(ii). We also get new and explicit models of $\DCF_{0,m}$ and $\CODF_m$: \begin{FACT}{Example}{ConstructionPowerSeriesCont} We continue example \ref{ConstructionPowerSeries}\ref{ConstructionPowerSeriesA}. If $K$ is a differential field, then the algebraic closure of the differentially large field $L=\bigcup _{n\in\N}K((\bt_1))\ldots((\bt_n))$ from \ref{ConstructionPowerSeries}\ref{ConstructionPowerSeriesA} is differentially closed. If $K$ is an ordered field and the order is extended to $L$ in some way, then the real closure of $L$ is a model of $\CODF_m$. \par Observe that these models are different from those obtained using iterated Puiseux series or generalized power series constructions in \ref{ConstructionPowerSeries}(ii). \end{FACT} \begin{FACT}{:Kolchin-Denseness of Rational Points in Differential Algebraic Groups.}{KolDenseInDiffGroups} \noindent In the classical case of a connected linear algebraic group $G$ over any field $F$ of characteristic zero, the Unirationality Theorem implies that the $F$-rational points of $G$ are Zariski-dense. In the differential situation the corresponding statement does not hold. For example, the linear differential algebraic group defined by $\delta x=x$ does in general not have a Kolchin dense set of rational points. However, in differentially large fields this is true again: \end{FACT} \begin{FACT}{Proposition}{groups} Assume $K$ is differentially large. If $G$ is a connected differential algebraic group over $K$, then the set of $K$-rational points of $G$, denoted $G(K)$, is Kolchin-dense in $G=G(\U)$. \end{FACT} \begin{proof} We verify \ref{CharDlargeI}\ref{Ivarieties}, hence it suffices to show that for infinitely many values of $r$ the jet $\jet_r G$ has a smooth $K$-rational point. By \cite[Corollary 4.2(ii)]{Pillay1997}, $G$ embeds over $K$ into a connected algebraic group $H$ defined over $K$. As we saw in \ref{diffalggeometry}, for each $r$, $\nabla_r G$ is a differential algebraic subgroup of $\tau_r H$. As a result, $\jet_r G$ is an algebraic subgroup of $\tau_r H$, and so $\jet_r G$ is smooth. If $e$ denotes the identity of $G$, which is a $K$-point, then, for each $r$, the $K$-point $\nabla_r(e)$ is a smooth point of $\jet_rG$. \end{proof} \begin{FACT}{:Pseudo differentially closed fields.}{psfDCF} Recall that a field $K$ is pseudo algebraically closed (PAC) if every absolutely irreducible algebraic variety over $K$ has a $K$-point. It is easy to see and well known that PAC fields are large and that the PAC property is equivalent to saying that $K$ is \ec\ in every regular field extension $L$ (meaning that $K$ is algebraically closed in $L$). From model theoretic literature one can formulate several notions of pseudo differentially closed fields. We show that they are all equivalent to the property ``PAC + differentially large". \par Let $K$ be a differential field. The following are equivalent. \begin{enumerate}[(i),itemsep=1ex] \item $K$ is PAC (as a field) and $K$ is differentially large. \item Every absolutely irreducible differential variety over $K$ has a differential $K$-point. Recall that a differential variety $V$ over $K$ is absolutely irreducible if it is irreducible in the Kolchin topology of a differential closure $K^{diff}$ and this is equivalent to saying that $V$ is irreducible over $K^{alg}$. \item $K$ is \ec\ in every differential field extension $L$ in which $K$ is R-\textit{regular} (i.e. $tp(a/K)$ is stationary for every tuple $a$ from $L$, where the type $tp(a/K)$ is with respect to the stable theory $DCF_{0,m}$). \item $K$ is \ec\ in every differential field extension $L$ in which $K$ is H-\textit{regular} (i.e. $K^{alg}\cap L=K$). \end{enumerate} \par \noindent If these conditions hold we call $K$ \notion{pseudo differentially closed}. \end{FACT} \begin{proof} We use \cite[Lemma 3.35]{Hoffma2019}, which says in our situation that a tuple $a$ from $\mathcal U$ (the monster model of $DCF_{0,m}$), the type $tp(a/K)$ is stationary if and only if the differential field extension $K\langle a\rangle$ over $K$ is $H$-regular. Clearly this characterization implies that H-regularity and R-regularity are equivalent. In particular, (iii) is equivalent to (iv). \par \smallskip\noindent (i)$\Ra $(iv) Let $K$ be PAC and differentially large. Let $L/K$ be an H-regular extension of $K$. Then $K$ is algebraically closed in $L$ as a field and because $K$ is PAC, it is \ec\ in $L$ as a field. Since $K$ is differentially large, it follows that it is \ec\ in $L$ as a differential field. \par \smallskip\noindent (iv)$\Ra $(ii) follows from the fact that a type $tp(a/K)$ is stationary if and only if the Kolchin-locus of $a$ over $K$ is absolutely irreducible: Let $V$ be an absolutely irreducible differential variety over $K$. Then the generic type $p=tp(a/K)$ of $V$ over $K$ is stationary, and hence, by the quoted characterization of stationarity, the differential field $L=K\langle a\rangle$ is an $H$-regular extension of $K$. By (iv), $K$ is \ec\ in $L$ as a differential field and so there is a differential $K$-point in $V$, as required. \par \smallskip\noindent (ii)$\Ra $(i) Suppose $V$ is a $K$-irreducible differential variety such that all jets of $V$ have a smooth $K$-point, then all these jets are absolutely Zariski irreducible (as they are Zariski $K$-irreducible and contain a smooth $K$-point). It follows that $V$ is absolutely irreducible. Hence, $V$ has a $K$-point. In fact, $V$ has Kolchin-dense many $K$-points; indeed, we can take any open differential subvariety $O$ of $V$ and argue similarly (using the fact that $K$ is large, as it is PAC, to produce smooth $K$-points in the jets of $O$). This shows (i) using the equivalence \ref{IDlarge}$\iff$\ref{Ivarieties} of \ref{CharDlargeI}. \end{proof} \par \iflongversion\LongStart \par \textcolor{red}{Next Corollary: Link to \cite{PilPol2006} and show how we get some (main?) results from there: Axiomatisation and Elementary equivalence} \begin{FACT}{Corollary}{} The theory of pseudo differentially closed fields is first order axiomatisable. If $M,L$ are pseudo-differentially closed fields and $K$ is a common differential subfield, then $M$ is elementary equivalent to $L$ over $K$ as differential fields just if they are elementary equivalent to $L$ over $K$ as fields. \end{FACT} \begin{proof} \underconstruction{} \begin{FACT}{:The elementary theory of a differentially large field}{} Let $M_1,M_2$ be differentially large fields and let $K$ be a common differential subfield. Suppose $M_1\equiv_KM_2$ as fields. Then $M_1\equiv_KM_2$ as differential fields. (Notice that in \cite[Theorem 6.2]{Tressl2005} this implication was shown with $\equiv_K$ replaced by $\equiv_{\forall,K}$.) \end{FACT} \begin{proof} We first show that for differentially large fields $K\subseteq L$ with $K_\alg\prec L_\alg$ we have $K\prec L$. Let $\SL$ be the language of rings and let $T_\alg$ be the $\SL$-theory of $K_\alg$. For each $\SL$-formula $\phi(x_1,\ldots,x_n)$ let $R_\phi$ be a new $n$-ary relation symbol, let $\SL^+$ be the extension of $\SL$ by all the $R_\phi$ and let $T_\alg^+$ be the $\SL^+$-theory extending $T$ by all the sentences $\forall x_1\ldots x_n(R_\phi \leftrightarrow \phi)$. Obviously $T_\alg^+$ has quantifier elimination in the language $\SL^+$ We claim that for every universal $\SL({\bar \delta})$-formula $\phi$, there is an existential $\SL^+({\bar \delta})$-formula $\psi$ such that $T_\alg^+\cup DL\vdash \phi\leftrightarrow \psi$. \par \medskip\centerline{\vbox{\hrule width 10em height 3pt}}\medskip \par By passing to elementary extensions of $M_1,M_2$ respectively we may assume that both are $|K|^+$-saturated. We will define a Back\ \&\ Forth system between $M_1$ and $M_2$ over $K$: Let \begin{align*} \DS=\{f:L_1\lra L_2\st &K\subseteq L_i\subseteq M_i\text{ intermediate differential field with }|L_i|=|K|\cr &\text{ and }f\text{ a differential }K\text{-algebra isomorphism, and}\cr &\tp^{M_{1,alg}}(L_1/K)=\tp^{M_{2,alg}}(L_2/K)\} \end{align*} \par \underconstruction{Not clear how to maintain $\tp^{M_{1,alg}}(L_1/K)=\tp^{M_{2,alg}}(L_2/K)$ after adding an element; the other conditions seem to be implementable using the twisted Taylor morphism starting with an extension $g$ of $f:L_1\lra L_2$ that is a field embedding (which exists because of $\tp^{M_{1,alg}}(L_1/K)=\tp^{M_{2,alg}}(L_2/K)$) } \section{Algebraic-Geometric Axiomatization of Large Differential Fields}\label{axiomsuc} \noindent One of our applications of the differential Weil descent constructed in Section~\ref{differentialweil} is to show that every algebraic extension of a large field that is a model of $\uc_m$ is again a model of $\uc_m$. This will be achieved in a similar manner to the fact that algebraic extensions of large fields are again large \cite[Proposition 2.1]{Pop1996}. The latter uses classical Weil descent and, in particular, the fact that the Weil functor preserves smoothness \cite[Appendix 2]{Oester1984}. Before presenting the applications we will review the theory $\uc_m$ and establish some useful characterizations. But first some preliminaries. \begin{FACT}{:Some preliminaries and conventions.}{} We fix a distinguished set of commuting derivations $\Delta=\{\delta_1,\dots,\delta_m\}$. We assume that all our fields are of characteristic zero. We work inside a large (saturated) differentially closed field $(\U,\Delta)$, and $K$ denotes a differential subfield of $\U$. A \notion{Kolchin-closed} subset of $\U^n$ is the common zero set of a set of differential polynomials over $\U$ in $n$ differential variables; such sets are also called \notion{affine differential varieties}. If the definining polynomials can be chosen with coefficients in $K$ we will say the set is \notion[]{defined over $K$}. The Kolchin-closed sets (defined over $K$) are the closed sets of a topology, called the \notion{Kolchin-topology} of $\U^n$ (over $K$). \par By a \notion{differential variety} $V$ we mean a topological space which has as finite open cover $V_1,\dots,V_s$ with each $V_i$ homeomorphic to an affine differential variety (inside some power of $\U$) such that the transition maps are regular as differential morphisms; see \cite[Chap. 1, section 7]{LeoSan2013}. We will say that the differential variety is over $K$ when all objects and morphisms can be defined over $K$. This definition also applies to our use of algebraic varieties, replacing Kolchin-closed with Zariski-closed in powers of $\U$ (recall that $\U$ is algebraically closed and a universal domain for algebraic geometry in Weil's ``foundations'' sense). \end{FACT} \begin{FACT}{Remark}{} Suppose $L/K$ is a finite field extension. Recall that the derivations $\delta_1,\dots,\delta_m$ extend uniquely from $K$ to $L$. \begin{enumerate} \item[(i)] Given a differential $L$-algebra $D$, by \ref{DiffWeilDescent}, there is a natural one-to-one correspondence between the differential $L$-points of $D$ and the differential $K$-points of $W^{\operatorname{diff}}(D)$. \item[(ii)] In the case when $D$ is the differential coordinate ring of an affine differential variety, say $D=L\{V\}$, and when $L$ has a $K$-basis $b_1,\dots,b_\ell$ of constants (meaning that $\delta_i(b_j)=0$ for all $i,j$), then a construction of the differential Weil descent $W^{\operatorname{diff}}(L\{V\})$ appears in \cite[\S5]{LeSMos2016}. However, a basis of constants does not always exist, see Example \ref{example1} below. \end{enumerate} \end{FACT} \begin{FACT}{Example}{example1} We work in the ordinary case $\Delta=\{\delta\}$. Let $K=\mathbb Q(t)$ with $\delta t=1$ and consider the finite extension $L=K(b)$ where $b^2=t$. Then the (unique) induced derivation on $L$ is given by $\delta b=\frac{1}{2b}=\frac{b}{2t}$. Fix the basis $\{1,b\}$ of $L$ as a $K$-module. Consider the differential variety $V$ given by $\delta x=0$ (i.e., $V$ is simply the constants of $\U$) viewed as a differential variety over $L$. The differential Weil descent $W^{\operatorname{diff}}(V)$ is obtained as follows; write $x$ as $x_1+x_2b$ and compute $$\delta(x_1+bx_2)=\delta x_1+(\delta b) x_2+b \delta x_2=\delta x_1+\frac{b}{2t} x_2+b \delta x_2=\delta x_1+\left(\frac{x_2}{2t}+\delta x_2\right)b.$$ Thus, $W^{\operatorname{diff}}(V)$ is the affine differential variety over $K$ given by the equations $$\delta x_1=0\quad \text{ and }\quad \delta x_2+\frac{x_2}{2t}=0.$$ Note that this is not contained in a product of the constants, as one might expect. Of course, if $\delta(b)$ were zero we would instead obtain the equations $\delta x_1=0$ and $\delta x_2=0$ (which would occur if $\delta$ were trivial on $K$, for instance). \end{FACT} \begin{FACT}{void}{} We fix integers $n>0$ and $r\geq 0$, and set $$ \Gamma_n(r) = \{(\xi,i) \in \mathbb N^m\times\{1,\dots,n\} \st \sum_{i=1}^m \xi_i \leq r\}. $$ \par \noindent We make use of prolongation spaces and recall the definition and some properties. The \notion[]{$r$-th nabla map} $\nabla_r:\U^n\to \U^{\alpha(n,r)}$ with $\alpha(n,r):=|\Gamma_n(r)|=n\cdot\binom{r+m}{m}$ is defined by $$\nabla_r(x)= (\delta^\xi x_i:\,(\xi,i)\in \Gamma_n(r)),$$ where $x=(x_1,\dots,x_n)$ and $\delta^\xi=\delta_1^{\xi_1}\cdots\delta_m^{\xi_m}$. We order the elements of the tuple $(\delta^\xi x_i:\,(\xi,i)\in \Gamma_n(r))$ according to the canonical orderly ranking of the indeterminates $\delta^\xi x_i$; that is, \begin{equation}\label{ordercanonical} \delta^{\xi}x_i< \delta^{\zeta}x_j \iff \left(\sum \xi_k,i,\xi_1,\dots,\xi_m\right)<\left(\sum \zeta_k,j,\zeta_1,\dots,\zeta_m\right) \end{equation} where the ordering on the right-hand-side is the lexicographic one. \par Let $\U_r:=\U[\epsilon_1,\dots,\epsilon_m]/(\epsilon_1,\dots,\epsilon_m)^{r+1}$ where the $\epsilon_i$'s are indeterminates, and let $e:\U\to \U_r$ denote the ring homomorphism $$x\mapsto \sum_{\xi\in\Gamma_1(r)}\frac{1}{\xi_1!\cdots\xi_m!}\; \delta^\xi(x)\; \epsilon_1^{\xi_1}\cdots\epsilon_m^{\xi_m}.$$ We call $e$ the exponential $\U$-algebra structure of $\U_r$. To distinguish between the standard and the exponential algebra structure on $\U_r$, we denote the latter by $\U_r^e$. \end{FACT} \begin{FACT}{Definition}{} Given an algebraic variety $X$ the $r$-th \notion{prolongation} $\tau X$ is the algebraic variety given by the taking the $\U$-rational points of the Weil descent of $X\times_\U \U_r^e$ from $\U_r$ to $\U$. Note that the base change $V\times_\U \U_r^e$ is with respect to the exponential structure while the Weil descent is with respect to the standard $\U$-algebra structure. \end{FACT} For details and properties of prolongation spaces we refer to \cite[\S 2]{MoPiSc2008}; for a more general presentation see \cite{MooSca2010}. In particular, it is pointed out there that the prolongation $\tau_r X$ always exist when $X$ is quasi-projective (an assumption that we will adhere to later on). A characterising feature of the prolongation is that for each point $a\in X=X(\U)$ we have $\nabla_r(a)\in \tau_r X$. Thus, the map $\nabla_{r}:X\to \tau_r X$ is a differential regular section of $\pi_r:\tau_r X\to X$ the canonical projection induced from the residue map $\U_r\to \U$. We note that if $X$ is defined over the differential field $K$ then $\tau_r X$ is defined over $K$ as well. \par In fact, $\tau_r$ as defined above is a functor from the category of algebraic varieties over $K$ to itself, and the maps $\pi_r:\tau_r X\to X$ and $\nabla_r:X\to \tau_r X$ are natural. The latter means that for any morphism of algebraic varieties $f:X\to Y$ we get \begin{equation}\label{natural} f\circ\pi_{r,X}=\pi_{r,Y}\circ\tau_r f \quad \text{ and }\quad \tau_r f\circ\nabla_{r,X}=\nabla_{r,Y}\circ f. \end{equation} If $G$ is an algebraic group, then $\tau_r G$ also has the structure of an algebraic group. Indeed, since $\tau_r$ commutes with products, the group structure is given by $$\tau_r(*):\tau_r G\times\tau_r G\to \tau_r G$$ where $*$ denotes multiplication in $G$. Moreover, by the right-most equality in \eqref{natural}, the map $\nabla_r:G\to\tau_r G$ is an injective group homomorphism. Hence, $\nabla_r(G)$ is a differential algebraic subgroup of $\tau_r G$. \par Assume that $V$ is a differential variety which is given as a differential subvariety of an algebraic variety $X$. We define the $r$-th jet of $V$ to be the Zariski-closure of the image of $V$ under the $r$-th nabla map $\nabla_r:X\to \tau_r X$; that is, $$\jet_r V=\overline{\nabla_r(V)}^{\operatorname{Zar}}\subseteq \tau_r X.$$ The jet sequence of $V$ is defined as $(\jet_r V:r\geq 0)$. Note that this sequence determines $V$, indeed $$V=\{a\in X: \nabla_r(a)\in \jet_r V \text{ for all $r\geq 0$}\}.$$ In the case when $V$ is affine, say a Kolchin-closed subset of $\U^n$, and defined by differential polynomials of order at most $r$, then $$V=\{a\in \U^n: \nabla_r(a)\in \jet_r V\}.$$ \begin{FACT}{:Assumption.}{} From now on we assume, whenever necessary for the existence of jets, that our differential varieties are given as differential subvarieties of quasi-projective algebraic varieties. Of course, in the affine case this is always the case. It is worth noting that for connected differential algebraic groups this is also true. Indeed, by \cite[Corollary 4.2(ii)]{Pillay1997} every such group embeds into a connected algebraic group and the latter is quasi-projective by Chevalley's theorem. \end{FACT} \medskip \begin{FACT}{:Reminder on characteristic sets.}{} We recall some of the theory of characteristic set of prime differential ideals of the differential polynomial ring $K\{x\}$ with $x=(x_1,\dots,x_n)$. For a detailed reference we refer the reader to \cite[Chapters I and IV]{Kolchi1973}. Let $f\in K\{x\}$ be nonconstant. The leader of $f$, denoted $v_f$, is the highest ranking algebraic indeterminate that appears in $f$ (according to the canonical orderly ranking of the indeterminates $\delta^\xi x_i$, as in the equivalence \eqref{ordercanonical}). The leading degree of $f$, denoted $d_f$, is the degree of $v_f$ in $f$. The rank of $f$, denoted $\operatorname{rk}(f)$, is the pair $(v_f,d_f)$. The set of ranks is ordered lexicographically. The separant of $f$, denoted $S_f$, is the formal partial derivative of $f$ with respect to $v_f$. The initial of $f$, denoted $I_f$, is the leading coefficient of $f$ when viewed as a polynomial in $v_f$. Note that both $S_f$ and $I_f$ have lower rank than $f$. Given a finite subset $\Lambda\subseteq K\{x\}\setminus K$, we set $H_\Lambda:=\prod_{f\in \Lambda}I_fS_f$. \par One says that $g\in K\{x\}$ is weakly reduced with respect to $f\in K\{x\}$ if no proper derivative of $v_f$ appears in $g$; if in addition the degree of $v_f$ in $g$ is strictly less than $d_f$ we say that $g$ is reduced with respect to $f$. A set $\Lambda\subseteq R\{x\}$ is said to be autoreduced if for any two distinct elements of $\Lambda$ one is reduced with respect to the other. Autoreduced sets are always finite, and we always write them in nondecreasing order by rank. The canonical orderly ranking on autoreduced sets is defined as follows: $\{g_1,\dots,g_r\}<\{f_1,\dots,f_s\}$ means that either there is $i\leq r,s$ such that $\operatorname{rk}(g_j)=\operatorname{rk}(f_j)$ for $j<i$ and $\operatorname{rk}(g_i)<\operatorname{rk}(f_i)$, or $r>s$ and $\operatorname{rk}(g_j)=\operatorname{rk}(f_j)$ for $j\leq s$. \par While it is not generally the case that prime differential ideals of $K\{x\}$ are finitely generated as differential ideals (though they are finitely generated as radical differential ideals), something close is true; they are determined by certain autoreduced subsets called characteristic sets. More precisely, if $P\subseteq K\{x\}$ is a prime differential ideal then a \notion{characteristic set} $\Lambda$ of $P$ is a minimal autoreduced subset of $P$ with respect to the ranking defined above. These minimal sets always exist, and determine the ideal $P$ in the sense that \begin{equation}\label{equchar} P=\{f\in K\{x\}:\; H_\Lambda^\ell f\in[\Lambda] \text{ for some $\ell\geq 0$} \}. \end{equation} The differential ideal on the right-hand-side is commonly denoted by $[\Lambda]:H_\Lambda^\infty$, where $[\Lambda]$ is the differential ideal generated by $\Lambda$ in $K\{x\}$. \end{FACT} \begin{FACT}{Fact}{reduced}\cite[Proposition 2.7]{Tressl2005} Suppose $\Lambda$ is a characteristic set of a prime differential ideal $P\subseteq K\{x\}$. If $f\neq 0$ is reduced with respect $\Lambda$, then $f$ is not in $P$. \end{FACT} \noindent Given $I\subseteq \U\{x\}$, let $\V(I)$ denote the zeroes (as differential solutions) of the elements of $I$ in $\U^n$. For a characteristic set $\Lambda$ of a prime differential ideal $P\subseteq K\{x\}$, the description \eqref{equchar} implies $$\V(P)\setminus \V(H_\Lambda)=\V(\Lambda)\setminus \V(H_\Lambda)$$ A consequence of Fact \ref{reduced} is that $H_\Lambda\notin P$, and hence the above equality says that $\V(P)$ and $\V(\Lambda)$ agree off a proper Kolchin-closed subset, namely $\V(H_\Lambda)$. \par \medskip \par We will need a bit more notation. We let $K\{x\}_{\leq r}$ denote the set of differential polynomials over $K$ of order at most $r$. On the other hand, letting $(x^{\xi}_i:(\xi,i)\in\N^m\times\{1,\dots,n\})$ be a collection of new variables, we set $$K\{x\}_{\leq r}^{\pol}=K[x^\xi _i:(\xi,i)\in \Gamma_n(r)].$$ More generally, if $\SS$ is a set of differential polynomials in $K\{x\}_{\leq r}$, we set $$\SS^{\pol}=\{f^{\pol}\in K\{x\}_{\leq r}^{\pol}:f\in \SS\}$$ where $f^{\pol}$ denotes the polynomial obtained by replacing the variables $\delta^\xi x_i$ in $f$ for the algebraic variables $x^\xi_i$. We also let $\V_r(\SS)$ denote the (algebraic) zero set of $\SS^{pol}$ in $\U^{\alpha(n,r)}$ where recall that $\alpha(n,r)=|\Gamma_n(r)|$. \begin{Remark}\label{onjets} If $V$ is an affine differential variety defined by the radical differential ideal $I\subseteq K\{x\}$, then for each $r$ the jet $\jet_r V$ has defining ideal given by $(I\cap K\{x\}_{\leq r})^{\pol}$. In other words, $$\jet_r V=\V_r(I\cap K\{x\}_{\leq r}).$$ \end{Remark} \par We can now recall the uniform companion theory $\UC_m$ of differential fields of characteristic zero with $m$ commuting derivations. For any set $\SS\subseteq K\{x\}$ we let $\SS^{(r)}$ denote the set of all $\delta^\xi f$ of order at most $r$ with $f\in \SS$. \begin{FACT}{Definition}{UCtheory}\cite{Tressl2005} A differential field $K$ is a model of $\UC_m$ if the following condition is satisfied: for every characteristic set $\Lambda$ of a prime differential ideal of $K\{x\}$, if $\V_r(\Lambda^{(r)})\setminus\V_r(H_\Lambda)\subseteq \U^{\alpha(n,r)}$ has a smooth $K$-point for some $r$ with $\Lambda\subseteq K\{x\}_{\leq r}$, then the differential variety $$\V(\Lambda)\setminus\V(H_\Lambda)\subseteq \U^{n}$$ has a (differential) $K$-point. \end{FACT} \begin{FACT}{Remark}{axiomrosenfeld} The fact that the class of differential fields that satisfy the above condition is first-order axiomatizable in the language of differential rings is the content of \cite[\S4]{Tressl2005}. The proof there relies heavily on Rosenfeld's Lemma which gives an algebraic characterization of characteristic sets of prime differential ideals \cite[Chapter IV, \S9]{Kolchi1973}. In Section \ref{axiomsuc} below we present an alternative (algebraic-geometric) axiomatization. \end{FACT} Next we prove two properties of characteristic sets of prime differential ideals that seem to be well known but to our knowledge are not explicitly stated elsewhere. \begin{FACT}{Lemma}{charsets1} Let $\Lambda$ be a characteristic set of a prime differential ideal $P\subseteq K\{x\}$. If $\Lambda\subseteq K\{x\}_{\leq r}$, then $$P\cap K\{x\}_{\leq r}=(\Lambda^{(r)}):H_\Lambda^\infty$$ where $(\Lambda^{(r)})$ denotes the ideal generated by $\Lambda^{(r)}$ in $K\{x\}_{\leq r}$. \end{FACT} \begin{proof} Since $P=[\Lambda]:H_\Lambda^\infty$, the containment $\supseteq$ is clear. Now let $f\in P\cap K\{x\}_{\leq r}$. Then $f$ is weakly reduced with respect to $\Lambda^{(r)}$. By the differential division algorithm, there is $g$ reduced with respect to $\Lambda^{(r)}$ and $\ell$ such that $$H_\Lambda^\ell \; f-g\in (\Lambda^{(r)}).$$ But then, as $f$ is in $P$, we get that $g\in [\Lambda]:H_\Lambda^\infty$. So, as $g$ is also reduced with respect to $\Lambda$, Fact~\ref{reduced} implies that $g=0$, and hence $f\in (\Lambda^{(r)}):H_\Lambda^\infty$. \end{proof} \begin{FACT}{Remark}{usefulfactjets} Suppose that $V$ is a $K$-irreducible affine differential variety with corresponding prime differential ideal $P\subseteq K\{x\}$. Putting together Remark~\ref{onjets} and Lemma~\ref{charsets1}, we obtain that if $\Lambda$ is a characteristic set of $P$ and $\Lambda\subseteq K\{x\}_{\leq r}$, then the defining ideal of $\Jet_r V$ in $K\{x\}_{\leq r}^{\pol}$ is $((\Lambda^{(r)}):H_\Lambda^{\infty})^{\pol}$. As a result, $$\Jet_r V\setminus \V_r(H_\Lambda)=\V_r(\Lambda^{(r)})\setminus\V_r(H_\Lambda).$$ \end{FACT} Recall that a field $F$ is \notion{existentially closed} (e.c. for short) in a field extension $L$ if every algebraic variety over $F$ with an $L$-point contains a $F$-point. The notion of e.c. for differential fields is defined similarly (in the category differential varieties). \begin{FACT}{Proposition}{charsets2} Suppose $\Lambda$ is a characteristic set of a prime differential ideal $P\subseteq K\{x\}$ and assume that $\Lambda\subseteq K\{x\}_{\leq r}$. Let $S=K\{x\}/P$ and $h=H_\Lambda/P\in S$, and let $$R=K\{x\}_{\leq r}/(\Lambda^{(r)}):H_\Lambda^\infty.$$ Then, $S_h$ is polynomial algebra over $R_h$. Consequently, $\operatorname{Frac}(R)$ is e.c. in $\operatorname{Frac}(S)$ as fields. \end{FACT} \begin{proof} Let $\Theta(x)^{>r}$ denote the set of derivatives $\delta^\xi x_i$ of order strictly larger than $r$. We thus have $$\Theta(x)^{>r}=\Theta_1(x)\cup \Theta_2(x)$$ where $\Theta_1(x)$ are elements of $\Theta^{>r}(x)$ that are not derivatives of any leader $v_f$ with $f\in \Lambda$, and $\Theta_2(x)=\Theta(x)^{>r}\setminus \Theta_1(x)$. We write $\bar \theta(x)$ for the coset of $\theta(x)$ in $S$. We claim that the elements of $\bar \Theta_1(x)\subseteq S$ are algebraically independent over $R_h$. Indeed, if there were $\theta_1(x),\dots,\theta_s(x)\in\Theta_1(x)$ such that $f(\bar\theta_1(x),\dots,\bar\theta_s(x))=0$ for some nonzero $f\in R_h[t_1,\dots,t_s]$, then for some $\ell$ we would get $H_\Lambda^\ell(x) f(\theta_1(x),\dots,\theta_s(x))\in P$. By the differential division algorithm, we can find $g$ reduced with respect to $\Lambda$ and $\ell'$ such that $$H_\Lambda^{\ell'}f(\theta_1(x),\dots,\theta_s(x))-g\in [\Lambda],$$ but since $\Lambda$ has order $\leq r$ and the $\theta_i(x)$'s are of order $>r$ and not a derivative a leader of $\Lambda$, we get that $g\neq 0$. But then $P$ would contain a nonzero element, namely $g$, that is reduced with respect to $\Lambda$, this contradicts Fact\ref{reduced}. \par We now prove that all the elements of $\bar\Theta_2(x)$ are in $R_h[\bar \Theta_1(x)]$. The result follows from this, as $S_h=R_h[\bar\Theta(x)]$. Let $\theta(x)\in \Theta_2(x)$. By the differential division algorithm, there is $g$ reduced with respect to $\Lambda$ and $\ell$ such that $H_\Lambda^\ell \theta(x)-g\in [\Lambda]$. But then, as $g\in K\{x\}_{\leq r}[\Theta_1(x)]$, we get $\bar\theta(x)\in R_h[\bar\Theta_1(x)]$. \end{proof} \par The above properties of characteristic sets are at the core of the proof of the Structure Theorem for differential algebras from \cite{Tressl2002}. In the next section we will use the following slightly different version of this theorem. \begin{FACT}{Theorem}{structuretheoremOLD} Let $B$ be a differential $K$-algebra that is differentially finitely generated and a domain. Then $B_h\cong_KA_h\otimes_K P$ where $A$ is a domain and a finitely generated $K$-algebra, $h\in A$, and $P$ is a polynomial algebra over $K$. \end{FACT} \begin{proof} By the assumptions, $B$ is of the form $K\{x\}/P$ for some tuple of differential indeterminates $x=(x_1,\dots,x_n)$ and $P$ a prime differential ideal of $K\{x\}$. Let $S$, $R$, $h$ and $\bar\Theta_1$ be as in the proof of Proposition~\ref{charsets2}, with $\Lambda$ a characteristic set of $P$, then if we set $A=R$ and $P=K[\bar\theta_1]$ the proposition yields that $S_h\cong_KA_h\otimes P$ with the desired properties. \end{proof} \par \textcolor{red}{Preliminaries above still need to be copied (and adapted) from Weil.tex} \par \LongEnd\else\fi \section{Algebraic-geometric axioms}\label{axiomsuc} \noindent In this last section we present algebraic-geometric axioms for differentially large fields in the spirit of the classical Pierce-Pillay axioms for differentially closed fields in one derivation \cite{PiePil1998} (cf. Remark~\ref{geoaxioms}(i) below). While this section might seem mostly of interest to model theorists, the general reader should keep in mind that Theorem~\ref{DiffLargeGeomAx} is a general statement on systems of algebraic PDEs that have solutions in differentially large fields. \par Our presentation here follows the recent algebraic-geometric axiomatization of differentially closed fields in several commuting derivations established in \cite{LeoSan2018}. In particular, we will use the recently developed theory of differential kernels for fields with several commuting derivations from~\cite{GusLeS2016}. One significant difference with the arguments in \cite{LeoSan2018} is that there one only requires the existence of regular realizations of differential kernels, while here we need the existence of principal realizations, see Remark~\ref{principalec} and Fact~\ref{useful}. We carry on the notation and conventions from previous sections. \par We use two different orders $\leq$ and $\unlhd$ on $\N^m \times \{1,\dots,n\}$. Given two elements $(\xi,i)$ and $(\tau,j)$ of $\N^m \times \{1,\dots,n\}$, we set $(\xi,i) \leq (\tau,j)$ if and only if $i = j$ and $\xi \leq \tau$ in the product order of $\N^m$. Furthermore, we set $(\xi,i) \unlhd (\tau,j)$ if and only if $$ (\sum \xi_k,i,\xi_1,\dots,\xi_{m}) \;\leq_{\text{lex}} \; (\sum \tau_k,j,\tau_1,\dots,\tau_m) $$ Note that if $x=(x_1,\ldots,x_{n})$ are differential indeterminates and we identify $(\xi,i)$ with $\delta^\xi x_i:=\delta_1^{\xi_1}\cdots\delta_m^{\xi_m}x_i$, then $\leq$ induces an order on the set of algebraic indeterminates given by $\delta^\xi x_i\leq \delta^\tau x_j$ if and only if $\delta^\tau x_j$ is a derivative of $\delta^\xi x_i$ (in particular this implies that $i=j$). On the other hand, the ordering $\unlhd$ induces the canonical orderly ranking on the set of algebraic indeterminates. \par We will look at field extensions of $K$ of the form \begin{equation}\label{extL} L:=K(a^\xi_i : (\xi,i) \in \Gamma_n(r)) \end{equation} for some fixed $r\geq 0$. Here we use $a^\xi_i$ as a way to index the generators of $L$ over $K$. The element $(\tau,j)\in \N^m \times \{1,\dots,n\}$ is said to be a leader of $L$ if there is $\eta\in \N^m$ with $\eta\leq \tau$ and $\sum \eta_k\leq r$ such that $a^\eta_j$ is algebraic over $K(a^\xi_i : (\xi,i) \lhd (\eta,j))$. A leader $(\tau,j)$ is a minimal leader of $L$ if there is no leader $(\xi,i)$ with $(\xi,i) < (\tau,j)$. Observe that the notions of leader and minimal leader make sense even when $r=\infty$. \par A (differential) kernel of length $r$ over $K$ is a field extension of the form $$L=K(a_i^{\xi}: (\xi,i)\in\Gamma_n(r))$$ such that there exist derivations $$D_k:K(a_i^\xi:(\xi,i)\in \Gamma_n(r-1))\to L$$ for $k=1,\dots,m$ extending $\delta_k$ and $D_ka_i^\xi=a_i^{\xi+{\bf k}}$ for all $(\xi,i)\in \Gamma_n(r-1)$, where ${\bf k}$ denotes the $m$-tuple whose $k$-th entry is one and zeroes elsewhere. \par Given a kernel $(L,D_1,\dots,D_k)$ of length $r$, we say that it has a \notion{prolongation of length $s\geq r$} if there is a kernel $(L',D_1',\dots,D_k')$ of length $s$ over $K$ such that $L'$ is a field extension of $L$ and each $D_k'$ extends $D_k$. We say that $(L,D_1,\dots,D_k)$ has a \notion{regular realization} if there is a differential field extension $(M,\Delta'=\{\delta_1',\dots,\delta_m'\})$ of $(K,\Delta=\{\delta_1,\dots,\delta_m\})$ such that $M$ is a field extension of $L$ and $\delta_k'a_i^\xi=a_i^{\xi+{\bf k}}$ for all $(\xi,i)\in \Gamma_n(r-1)$ and $k=1,\dots ,m$. In this case we say that $g:=(a_1^{\bf 0},\dots,a_n^{\bf 0})$ is a regular realization of $L$. If in addition the minimal leaders of $L$ and those of the differential field $K\langle g\rangle$ coincide we say that $g$ is a \notion{principal realization} of $L$. \begin{FACT}{Remark}{principalec} Note that if $g$ is a principal realization of the differential kernel $L$, then $L$ is existentially closed in $K\langle g \rangle$ as fields. Indeed, since the minimal leaders of $L$ and $K\langle g\rangle$ coincide, for every $(\xi,i)\in \N^m \times \{1,\dots,n\}$ we have that either $\delta^\xi g_i$ is in $L$ or it is algebraically independent from $K(\delta^\eta g_j:(\eta,j)\lhd (\xi,i))$. In other words, the differential ring generated by $g$ over $L$, namely $L\{ g\}$, is a polynomial ring over $L$. The claim follows. \end{FACT} \noindent In general, it is not the case that every kernel has a principal realization (not even regular). In \cite{GusLeS2016}, an upper bound $C_{r,m}^n$ was obtained for the length of a prolongation of a kernel that guarantees the existence of a principal realization. This bound depends only on the data $(r,m,n)$ and is constructed recursively as follows:\looseness=-1 $$C_{0,m}^1=0, \quad\; C_{r,m}^1=A(m-1,C_{r-1,m}^1), \quad \text{ and } \quad C_{r,m}^n=C_{C_{r,m}^{n-1},m}^1,$$ where $A(x,y)$ is the Ackermann function. For example, $$C_{r,1}^n=r, \quad\; C_{r,2}^n=2^n r \quad \text{ and }\quad C_{r,3}^1=3(2^r-1).$$ In \cite[Theorem 18]{GusLeS2016}, it is proved that \begin{FACT}{Fact}{useful} If a differential kernel $L=K(a_i^\xi:(\xi,i)\in\Gamma_n(r))$ of length $r$ has a prolongation of length $C_{r,m}^n$, then there is $r\leq h\leq C_{r,m}^n$ such that the differential kernel $K(a_i^\xi:(\xi,i)\in \Gamma_n(h))$ has a principal realization. \end{FACT} \begin{FACT}{Remark}{} Note that in the ordinary case $\Delta=\{\delta\}$ (i.e., $m=1$), we have $C_{r,1}^n=r$ by definition, and so the fact above shows that in this case every differential kernel has a principal realization (this is a classical result of Lando \cite{Lando1970}). \end{FACT} \noindent The fact above is the key to our algebraic-geometric axiomatization of differential largeness. We need some additional notation. For a given positive integer $n$, we set $$\alpha(n)=n\cdot\binom{C_{1,m}^n+m}{m}\quad\text{ and }\quad \beta(n)=n\cdot\binom{C_{1,m}^n-1+m}{m}.$$ We write $\pi:\U^{\alpha(n)}\to \U^{\beta(n)}$ for the projection onto the first $\beta(n)$ coordinates; i.e., setting $(x_i^{\xi})_{(\xi,i)\in\Gamma_n(C_{1,m}^n)}$ to be coordinates for $\U^{\alpha(n)}$ then $\pi$ is the map $$(x_i^{\xi})_{(\xi,i)\in\Gamma_n(C_{1,m}^n)}\mapsto (x_i^{\xi})_{(\xi,i)\in\Gamma(C_{1,m}^n-1)}.$$ It is worth noting here that $\alpha(n)=|\Gamma_n(C_{1,m}^n)|$ and $\beta(n)=|\Gamma_n(C_{1,m}^n-1)|$. We also use the projection $\psi:\U^{\alpha(n)}\to \U^{n\cdot(m+1)}$ onto the first $n\cdot(m+1)$ coordinates, that is,\looseness=-1 $$(x_i^{\xi})_{(\xi,i)\in\Gamma_n(C_{1,m}^n)}\mapsto (x_i^{\xi})_{(\xi,i)\in\Gamma_n(1)}.$$ Finally, we use the embedding $\phi:\U^{\alpha(n)}\to \U^{\beta(n)\cdot(m+1)}$ given by \begin{align*} (x_i^{\xi})_{(\xi,i)\in\Gamma_n(C_{1,m}^n)}\mapsto \biggl((x_i^{\xi})_{(\xi,i)\in\Gamma_n(C_{1,m}^n-1)},&(x_i^{\xi+{\bf 1}})_{(\xi,i)\in\Gamma_n(C_{1,m}^n-1)},\ldots\cr &\ldots,(x_i^{\xi+{\bf m}})_{(\xi,i)\in\Gamma_n(C_{1,m}^n-1)}\biggr). \end{align*} \par \noindent Recall from\ref{diffalggeometry} that, given a Zariski-constructible set $X$ of $\U^n$, the first-prolongation of $X$ is denoted by $\tau X=\tau_1 X\subseteq \U^{n(m+1)}$. For the first-prolongation it is easy to give the defining equations: $\tau(X)$ is the Zariski-constructible set given by the conditions $$x\in X, \quad \text{ and }\quad \sum_{i=1}^n\frac{\partial f_j}{\partial x_i}(x)\cdot y_{i,k} +f_j^{\delta_k}(x)=0\; \text{ for } 1\leq j\leq s, \; 1\leq k\leq m$$ where $f_1,\dots,f_s$ are generators of the ideal of polynomials over $\U$ vanishing at $X$, and each $f_j^{\delta_k}$ is obtained by applying $\delta_k$ to the coefficients of $f_j$. Note that $(a,\delta_1 a,\dots,\delta_m a)\in \tau X$ for all $a\in X$. Further, if $X$ is defined over the differential field $K$ then so is $\tau X$. \begin{FACT}{Theorem}{DiffLargeGeomAx} Assume $K$ is a differential field that is large as a field. Then, $K$ is differentially large if and only \begin{enumerate} \item[($\diamondsuit$)] for every $K$-irreducible Zariski-closed set $W$ of $\U^{\alpha(n)}$ with a smooth $K$-point such that $\phi(W)\subseteq \tau(\pi(W))$, the set of $K$-points of $\psi(W)$ of the form $(a,\delta_1 a,\dots,\delta_m a)$ is Zariski-dense in $\psi(W)$. \end{enumerate} \end{FACT} \begin{proof} The proof follows the strategy of \cite{LeoSan2018}, but here regular realizations are replaced by principal realizations with the appropriate adaptations. As the set up is technically somewhat intricate we give details. \par Assume $K$ is differentially large. Let $W$ be as in condition ($\diamondsuit$), we must find Zariski-dense many $K$-points in $\psi(W)$ of the form $(a,\delta_1 a,\dots,\delta_m a)$. Let $b=(b_i^{\xi})_{(\xi,i)\in\Gamma_n(C_{1,m}^n)}$ be a Zariski-generic point of $W$ over $K$. Then $(b_i^{\xi})_{(\xi,i)\in\Gamma_n(C_{1,m}^n-1)}$ is a Zariski-generic point of $\pi(W)$ over $K$, and \begin{align*} \phi(b)&=\left((b_i^{\xi})_{(\xi,i)\in\Gamma_n(C_{1,m}^n-1)},(b_i^{\xi+{\bf 1}})_{(\xi,i)\in\Gamma_n(C_{1,m}^n-1)},\dots,(b_i^{\xi+{\bf m}})_{(\xi,i)\in\Gamma_n(C_{1,m}^n-1)}\right)\cr &\in \tau (\pi(W)) \end{align*} \noindent By the standard argument for extending derivations (see \cite[Chapter~7, Theorem~5.1]{Lang2002}, for instance), there are derivations $$D'_k:K(b_i^\xi:(\xi,i)\in \Gamma_n(C_{1,m}^n-1))\to K(b_i^\xi:(\xi,i)\in\Gamma_n(C_{1,m}^n))$$ for $k=1,\dots,m$ extending $\delta_k$ and such that $D'_kb_i^\xi=b_i^{\xi+{\bf k}}$ for all $(\xi,i)\in\Gamma_n(C_{1,m}^n-1)$. Thus, $L'=K(b_i^\xi:(\xi,i)\in\Gamma_n(C_{1,m}^n))$ is a differential kernel over $K$ and, moreover, it is a prolongation of length $C_{1,m}^n$ of the differential kernel $L=K(b_i^\xi:(\xi,i)\in\Gamma_n(1))$ of length 1 with $D_k=D_k'|_{L}$. By Fact \ref{useful}, there is $r\leq h\leq C_{1,m}^n$ such that $L''=K(b_i^\xi:(\xi,i)\in\Gamma_n(h))$ has a principal realization; in particular, there is a differential field extension $(M,\Delta')$ of $(K,\Delta)$ containing $L''$ such that $\delta_k' b^{\bf 0}=b^{\bf k}$, where $b^{\bf 0}=(b_1^{\bf 0},\dots,b_n^{\bf 0})$ and similarly for $b^{\bf k}$. Then \[ (b^{\bf 0},\delta_1' b^{\bf 0}, \dots,\delta_m' b^{\bf 0}) \text{ is a generic point of } \psi(W) \text{ over } K\leqno{(*)} \] Now, since $W$ has a smooth $K$-point and $K$ is large, $K$ is e.c. in $L'$ as fields; in particular, $K$ is e.c. in $L''$ as fields. By Remark~\ref{principalec}, $L''$ is e.c. in the differential field $K\langle b^{\bf 0}\rangle$ as fields, and so $K$ is e.c. in $K\langle b^{\bf 0}\rangle$ as fields. Since $K$ is differentially large, the latter implies that $K$ is e.c. in $K\langle b^{\bf 0}\rangle$ as differential fields as well. The conclusion now follows using $(*)$. \par For the converse, assume $K$ is e.c. as field in a differential field extension $F$. We must show that $K$ is also e.c. in $F$ as differential field. Let $\rho(x)$ be a quantifier-free formula over $K$ (in the language of differential rings with $m$ derivations) in variables $x=(x_1,\dots,x_t)$ with a realization $c$ in $F$. We may write $$\rho(x)=\gamma(\delta^\xi x_i:(\xi,i)\in\Gamma_t(r))$$ where $\gamma((x^\xi)_{(\xi,i)\in\Gamma_t(r)})$ is a quantifier-free formula in the language of rings over $K$ for some $r$. If $r=0$, then $\rho$ is a formula in the language of rings, and so $\rho(x)$ has a realization in $K$ since $K$ is e.c. in $F$ as a field. Now assume $r>0$. Let $n:=t\cdot\binom{r-1+m}{m}$, $d:=(\delta^\xi c_i)_{(\xi,i)\in\Gamma_t(r-1)}$, and $$W:=\operatorname{Zar-loc}_K(\delta^\xi d_i:(\xi,i)\in\Gamma_n(C_{1,m}^n))\subseteq \U^{\alpha(n)}.$$ We have that $\phi(W)\subseteq \tau(\pi(W))$. Moreover, since $W$ has a smooth $F$-point (namely $(\delta^\xi d_i)_{(\xi,i)\in\Gamma_n(C_{1,m}^n)}$) and $K$ is e.c. in $F$ as fields, $W$ has a smooth $K$-point. By ($\diamondsuit$), there is $a=(a^\xi_i)_{(\xi,i)\in\Gamma_t(r-1)}\in K^n$ such that $(a,\delta_1 a,\dots,\delta_m a)\in \psi(W)$. This implies that $a^{\xi}_i=\delta^\xi a_i^{\bf 0}$ for all $(\xi,i)\in\Gamma_t(r-1)$. Thus, $$(\delta^\xi a_i^{\bf 0})_{(\xi,i)\in\Gamma_t(r)}\in\operatorname{Zar-loc}_K((\delta^\xi c_i)_{(\xi,i)\in\Gamma_t(r)})\subseteq \U^{t\cdot\binom{r+m}{m}},$$ and so, since $(\delta^\xi c_i)_{(\xi,i)\in\Gamma_t(r)}$ realizes $\gamma$, the point $(\delta^\xi a_i^{\bf 0})_{(\xi,i)\in\Gamma_t(r)}$ also realizes $\gamma$. Consequently, $K\models \rho(a^{\bf 0})$, as desired. \end{proof} \par \noindent In the ordinary case ($m=1$) we get the values $\alpha(n)=2n$ and $\beta(n)=n$. Also, in this case, $\pi:\U^{2n}\to \U^n$ is just the projection onto the first $n$ coordinates, and $\psi, \phi:\U^{2n}\to \U^{2n}$ are both the identity map. We thus get the following \begin{FACT}{Corollary}{ordinary} Assume that $(K,\delta)$ is an ordinary differential field of characteristic zero which is large as a field. Then, $(K,\delta)$ is differentially large if and only \begin{enumerate} \item[($\diamondsuit'$)] for every $K$-irreducible Zariski-closed set $W$ of $\U^{2n}$ with a smooth $K$-point such that $W\subseteq \tau_\delta(\pi(W))$, the set of $K$-points of $W$ of the form $(a,\delta a)$ is Zariski-dense in $W$. \end{enumerate} \end{FACT} \begin{FACT}{Remark}{geoaxioms} \begin{enumerate}[(i)] \item If $K$ is algebraically closed of characteristic 0, then Corollary~\ref{ordinary} yields the classical algebraic-geometric axiomatization of $\operatorname{DCF_0}$ given by Pierce and Pillay in \cite{PiePil1998}. \item If $K$ has a model complete theory $T$ in the language of fields and if $K$ is large, then Corollary~\ref{ordinary} yields a slight variation of the geometric axiomatization of $T_D$ given by Brouette, Cousins, Pillay and Point in \cite[Lemma 1.6]{BCPP2017}. \item For large and topological fields with a single derivation, an alternative description of differentially large fields with reference to the topology may be found in \cite{GuzRiv2006}. \end{enumerate} \end{FACT} \renewcommand\href[3][]{}  \par \end{document}